\documentclass[aos,MSNbibl,nameyear,seceqn,dvips]{arximspdf}

\doi{10.1214/13-AOS1106} %kopijuoti is PTS
\volume{41}
\issue{3}
\pubyear{2013}
\firstpage{1232}
\lastpage{1259}

\makeatletter

\newcommand{\eqref}[1]{(\ref{#1})}
\newtheorem{theorem}{Theorem}[section]
\newtheorem{lemma}[theorem]{Lemma}

\newproclaim{remark}{Remark}

\newcommand{\pr}{\mathbf{P}}
\newcommand{\E}{\mathbf{E}}
\newcommand{\Var}{\operatorname{Var}}
\newcommand{\ps}{\mathbf{P}_{*}}
\newcommand{\Es}{{\mathbf{E}_*}}

\newcommand{\alb}{\bolds{\alpha}}
\newcommand{\bb}{\bolds\beta}
\newcommand{\glb}{\bolds{\gamma}}
\newcommand{\xib}{\bolds{\xi}}
\newcommand{\etab}{\bolds{\eta}}
\newcommand{\bbln}{\widehat{\bolds{\beta}}_n}
\newcommand{\bbtn}{\widetilde{\bolds{\beta}}_n}

\newcommand{\bbf}{\mathbf{b}}
\newcommand{\xb}{\mathbf{x}}
\newcommand{\bfx}{\mathbf{x}}
\newcommand{\yb}{\mathbf{y}}
\newcommand{\sbf}{\mathbf{s}}
\newcommand{\ub}{\mathbf{u}}
\newcommand{\tb}{\mathbf{t}}

\newcommand{\A}{\mathbf{A}}
\newcommand{\C}{\mathbf{C}}
\newcommand{\D}{\mathbf{D}}
\newcommand{\Q}{\mathbf{Q}}
\newcommand{\R}{\mathbf{R}}

\newcommand{\T}{\mathbf{T}}
\newcommand{\V}{\mathbf{V}}
\newcommand{\W}{\mathbf{W}}

\newcommand{\cC}{\mathcal{C}}

\newcommand{\Zero}{\mathbf{0}}

\newcommand{\tr}{ \prime}

\newcommand{\al}{\alpha}

\newcommand{\ep}{\varepsilon} %%%%%%% epsilon is substituted as
\newcommand{\ga}{\gamma}
\newcommand{\La}{\Lambda}

\newcommand{\si}{\sigma}
\newcommand{\Up}{\Upsilon}

\newcommand{\nti}{n\rightarrow\infty}
\newcommand{\rai}{\rightarrow\infty}
\newcommand{\raw}{\rightarrow}
\makeatother

\begin{document}
\begin{frontmatter}

\title{Rates of convergence of the Adaptive LASSO estimators
to the Oracle distribution and higher order refinements by the bootstrap}
\runtitle{Edgeworth expansions for Adaptive LASSO estimators}

\begin{aug}
\author[A]{\fnms{A.} \snm{Chatterjee}\thanksref{t1}\corref{}\ead[label=e1]{cha@isid.ac.in}}
\and
\author[B]{\fnms{S.~N.} \snm{Lahiri}\thanksref{T1}\ead[label=e2]{snlahiri@ncsu.edu}}

\thankstext{t1}{Supported in part by the VI-MSS program of Department
of Science and Technology,
Government of India, and the Statistical and Applied Mathematical
Sciences Institute (SAMSI), NC, USA.}
\thankstext{T1}{Supported in part by NSF Grant DMS-10-07703
and NSA Grant H98230-11-1-0130.
On leave from Texas A\&M University.}
\runauthor{A. Chatterjee and S.~N. Lahiri}

\affiliation{Indian Statistical Institute and North Carolina State University}

\address[A]{Statistics and Mathematics Unit\\
Indian Statistical Institute\\
New Delhi 110067\\
India\\
\printead{e1}}
\address[B]{Department of Statistics\\
North Carolina State University\\
Raleigh, North Carolina 27695\\
USA\\
\printead{e2}}
\end{aug}

% HISTORY:
\received{\smonth{2} \syear{2012}}
\revised{\smonth{1} \syear{2013}}

% ABSTRACT
%
\begin{abstract}
Zou [\textit{J. Amer. Statist. Assoc.} \textbf{101} (2006) 1418--1429]
proposed the Adaptive LASSO
(ALASSO) method for simultaneous variable selection and
estimation of the regression parameters, and
established its oracle property. In this paper, we
investigate the rate of convergence of the
ALASSO estimator to the oracle distribution
when the dimension of the regression parameters
%is undounded.
may grow to infinity with the sample size.
It is shown that the rate critically depends on the choices of
the penalty parameter and the initial estimator,
among other factors, and that confidence intervals
(CIs) based on the oracle limit law often have poor
coverage accuracy. As an alternative, we consider
% higher order properties of
the residual bootstrap method for the ALASSO
estimators %of the nonzero regression parameters
that has been recently shown to be consistent; cf.
Chatterjee and Lahiri [\textit{J. Amer. Statist. Assoc.} \textbf{106}
(2011a) 608--625].
%;\textit{J. Amer. Statist. Assoc.}).
We show that the bootstrap applied
to a suitable studentized version of the
ALASSO estimator achieves second-order correctness,
even when the dimension of the regression parameters
is unbounded.
%Further, using higher order expansions, we
%develop a new bias-adjusted studentized version of
%the ALASSO estimator and show that the bootstrap
%approximation to the distribution of this modified
%pivot attains the same optimal level of accuracy as
%in the traditional case where the dimensions
% of the parameters id fixed
Results from a moderately large simulation study show
marked improvement in coverage accuracy for the bootstrap
CIs over the oracle based CIs.
\end{abstract}

% KEYWORDS


\end{frontmatter}

%s1 #&#
%s1 ###
\section{Introduction}
\label{sec1}

Consider the regression model
%
%e1.1 #&#
%
%e1.1 ###
\begin{equation}
y_i = \xb^{\tr}_i\bb+ \ep_i,\qquad
i=1,\ldots,n, \label{true-model}
\end{equation}
where $y_i$ is the response,
$\xb_i={ (x_{i,1},\ldots,x_{i,p} )}{}^{\prime}$
is a $p$ dimensional covariate vector,
$\bb={ (\beta_1,\ldots,\beta_p )}{}^{\prime}$
is the regression parameter and $ \{\ep_i\dvtx i=1,\ldots,n \}$
are independent and identically\vadjust{\goodbreak} distributed (i.i.d.) errors.
Let $\bbtn$ denote a root-$n$ consistent estimator of $\bb$, such as
the ordinary least squares (OLS) estimator of $\bb$. The Adaptive Lasso
(ALASSO) estimator of $\bb$ is defined as the minimizer of
the weighted $\ell_1$-penalized least squares criterion function,
%
%e1.2 #&#
%
%e1.2 ###
\begin{equation}
\bbln= \mathop{\operatorname{argmin}}_{\ub\in\mathbb{R}^p}\sum
_{i=1}^n { \bigl(y_i-
\xb^{\prime}_i\ub \bigr)}^2 +\lambda_n
\sum_{j=1}^p \frac{|u_j|}{{|\tilde{\beta}_{j,n}|}^\gamma},
\label{alasso-def}
\end{equation}
where $\lambda_n>0$ is a regularization parameter, $\gamma>0$
and $\tilde{\beta}_{j,n}$ is the $j$th component of $\bbtn$. The
ALASSO provides an improvement over the LASSO and related
bridge estimators that often require strong regularity
conditions on the design vectors $\xb_i$'s for consistent
variable selection and that have nontrivial bias in the
selected nonzero components; cf. \citet{knight-fu}, \citet{fan-li-01},
\citet{yuan-lin-07}, \citet{zhao-yu-06}. To
highlight some of the key properties
of the ALASSO, suppose
for the time being, that the first $p_0$ components of the true
regression parameter $\bb$ are nonzero and the last $(p-p_0)$
components are zero, where $1\leq p_0<p$. Let
$\widetilde{\mathcal{I}}_n=\{j\dvtx1\leq j\leq p,
\widehat{\beta}_{j,n}\neq0\}$ denote the variables
selected by the ALASSO, where $\widehat{\beta}_{j,n}$ is
the $j$th component of $\bbln$. \citet{zou-06} showed
that under some mild regularity conditions,
for fixed $p$, as $n\raw\infty$,
%
%e1.3 #&#
%
%e1.3 ###
\begin{equation}
\label{c-v} \pr(\widetilde{\mathcal{I}}_n =\mathcal{I}_n
) \raw1\quad \mbox{and}\quad \sqrt{n} \bigl(\widehat{\bb}^{(1)}_n-
\bb^{(1)} \bigr)\stackrel{d} {\rightarrow} N \bigl(\Zero,
\sigma^2\C^{-1}_{11} \bigr),
\end{equation}
%
%Thus, the ALASSO method has the `consistent variable selection'
%property. \citet{zou-06} also showed that
% the asymptotic distribution of the ALASSO estimators
%of the selected variables
% satisfies:
%N\left(\Zero,\sigma^2\C^{-1}_{11}\right),
where
$ \mathcal{I}_n
= \{1,\ldots,p_0\}$,
$\widehat{\bb}^{(1)}_n= (\widehat{\beta}_{1,n},
\ldots,\widehat{\beta}_{p_0,n} )$,
$\bb^{(1)}= (\beta_1,\ldots,\beta_{p_0} )$ and $\C
_{11}$ is
the upper left $p_0\times p_0$ submatrix of $\C\equiv
\lim_{\nti}n^{-1}\sum_{i=1}^n \xb_i\xb^{\prime}_i$.
% Note that if
% we had {\it{apriori}} knowledge that $\beta_{p_0+1}=\cdots=
%of the OLS estimator of $\bb^{(1)}$ coincides with
%the limit distribution in \eqref{c-v}.
Thus, the ALASSO method enjoys
the \textit{oracle property} [cf. \citet{fan-li-01}], that is,
it can correctly identify the set of nonzero components of $\bb$,
with probability tending to $1$
and at the same time, estimate the nonzero components accurately,
with the same precision as that of the OLS method, in the limit.

%Asymptotic normality of the ALASSO estimators of the nonzero
%components is an important first step as it makes it possible to
% do statistical inference on the nonzero
% regression parameters, following variable selection.
%However,
Although the oracle property of the ALASSO estimators
allows one to carry out statistical inference on the nonzero
regression parameters, following variable selection,
accuracy of
%the approximation by the limiting normal
% distribution is an important issue, as it determines the quality
of the resulting inference remains unknown.
In this paper, we investigate the rate of convergence of
$\sqrt{n}(\widehat{\bb}^{(1)}_n-\bb^{(1)})$ to the oracle
limit and show
that the penalization term in \eqref{alasso-def}
%of the ALASSO estimator
induces a substantial amount of bias
which, although vanishes asymptotically, can lead to a poor
rate of convergence.
% to the normal limit.
As a result,
large sample inference based on the oracle distribution
is not very accurate. As an alternative, we consider
the bootstrap method
or more precisely,
the residual bootstrap method [cf.
\citet{efron-aos-79},
\citet{freedman-81}], that is,
the most common version of the bootstrap in a regression model
like \eqref{true-model}. Recently,
Chatterjee and
Lahiri (\citeyear{ac-snl-Lasso-08,ac-snl-jasa-11}) showed that while
the residual bootstrap drastically fails for the LASSO.
Rather
surprisingly,
it provides a valid approximation to the distribution of
the centered and scaled ALASSO-estimator.
% under some mild conditions.
Notwithstanding its success in
capturing the first order limit, the accuracy of the bootstrap
for the ALASSO remains unknown. In this paper,
we also study the rate of bootstrap approximation to the
distribution of the ALASSO estimators, with and without
studentization, and develop ways to improve it,
all in the more general framework where the
number of regression parameters $p=p_n$
is allowed to go to infinity with the sample
size $n$.

To describe the main findings of the paper, consider
\eqref{true-model} where $p$, $\xb_i$'s and~$\bb$ are allowed
to depend on $n$ (but we often suppress the subscript $n$
to ease notation) and let
$
\T_n = \sqrt{n}\D_n(\bbln-\bb)
$,
where $\D_n$ is a known $q\times p$ matrix with
$\operatorname{tr}(\D_n\D^{\prime}_n)=O(1)$ and $q\in\mathbb{N}=\{
1,2,\ldots\}$
is an integer, {\textit{not}} depending on $n$. Thus, $\T_n$ is the vector
of $q$ linear functions of
$n^{1/2}(\bbln-\bb)$. Under the regularity conditions of
Section~\ref{sec3}, $\{\T_n\dvtx n\geq1\}$ is asymptotically
normal with mean zero and $q\times q$ asymptotic variance
${\bolds{\Sigma}}_n$ (say). We
consider
the error of oracle-based normal approximation,
\[
\Delta_n  \equiv\sup_{B\in\cC_q} \bigl|\pr(\T_n
\in B) -\Phi(B;\bolds{\Sigma}_n) \bigr|,
\]
where,
for $k\geq1$,
$\cC_k$ is the collection of all convex measurable
subsets of $\mathbb{R}^k$ and
$\Phi(\cdot;\A)$ is the Gaussian measure on
$\mathbb{R}^k$ with mean zero and $k\times k$ covariance
matrix~$\A$. Theorem~\ref{thm3-1} below gives an upper bound
on $\Delta_n$,
%on the rate of convergence of the difference:
%
%e1.4 #&#
%
%e1.4 ###
\begin{equation}
\Delta_n \leq\operatorname{const} \bigl[n^{-1/2}+\|
\bbf_n\| +c_n \bigr], \label{o-rate}
\end{equation}
where $\bbf_n$ is a bias term that results from the
penalization scheme in \eqref{alasso-def} and where
$c_n\in(0,\infty)$ is determined by
the initial $\sqrt{n}$-consistent estimator $\bbtn$
and the tuning parameter $\gamma$ in
\eqref{alasso-def}. The magnitude of both these terms
critically depend on the choice of the penalization parameter
$\lambda_n$ and the exponent $\gamma$, and either of them can
make the error rate sub-optimal, that is, worse than the
rate $O(n^{-1/2})$ that is attained by the
oracle based OLS estimator. Further, Theorem~\ref{thm3-1b}
shows that under some additional mild conditions, the
rate in \eqref{o-rate} is \textit{exact}, that is,
$\Delta_n$ is also \textit{bounded below} by
a constant multiple of the sum of the
three terms on the right-hand side of \eqref{o-rate}.
Therefore, it follows that although the ALASSO estimator converges
to the oracle distribution in the limit, the convergence rate
can be sub-optimal.
% and in general, this cannot be improved.
A direct implication of this result is that large sample tests
and CIs based on the normal limit law of the
ALASSO estimator
%of the nonzero parameters
may perform poorly, depending on the choice of the regularization
parameters $\lambda_n$ and $\gamma$. The simulation results
of Section~\ref{sec5} confirm this finite samples.

Next we consider properties of bootstrap approximations to the
distributions of
$\T_n$ and $\R_n$, a \textit{computationally simple} studentized
version of $\T_n$, given by
$
\R_n = \frac{\T_n}{\widehat{\sigma}_n}
$,
where $\widehat{\sigma}^2_n$ is the sample variance of the
ALASSO based residuals. Here we use a scalar studentizing factor
instead of the usual matrix factor [cf. \citet{lahiri-1994-sankhya}]
to reduce the computational burden. Fortunately, this does not
impact the accuracy of the bootstrap approximation as $\sigma^2$ is
the only unknown population parameter in the limit
distribution of $\T_n$. Theorem~\ref{thm3-2} below shows
that under fairly general conditions,
the rate of bootstrap approximation to the distribution of
$\T_n$ is $O_p(n^{-1/2})$.
Thus, the bootstrap corrects for the effects of
$\|\bbf_n\|$ and $c_n$ in \eqref{o-rate}, and produces a more
``accurate'' approximation to the distribution of $\T_n$
than the oracle based normal approximation. % in \eqref{o-rate}.
As a consequence, bootstrap percentile CIs based on
the ALASSO have a better performance compared to the large sample
normal CIs based on the oracle.

The results on the studentized statistic $\R_n$ are more encouraging.
Theorem~\ref{thm-st} shows that the bootstrap applied to $\R_n$ has
an error rate of $o_p(n^{-1/2})$ which outperforms the best possible
rate, namely $O(n^{-1/2})$ of normal approximation, irrespective of
the order of the terms $\|\bbf_n\|$ and $c_n$ in \eqref{o-rate}.
%The proof of Theorem~\ref{thm-st} shows that the (conditional)
%distribution of the bootstrapped statistic is able to capture the
% smaller order terms corresponding to $n^{-1/2}$, $\|\bbf_n\|$ and
%$c_n$,
% and as a consequence,
Thus, the bootstrap applied to the studentized
statistic $\R_n$ achieves \textit{second order correctness}. In contrast,
the normal approximation to the distribution of $\R_n$
%fails to
%account for each of these three terms and
has an error of
the order $O(n^{-1/2}+\|\bbf_n\|+c_n)$, as in the case of~$\T_n$.
As a result, bootstrap percentile-$t$ CIs based on $\R_n$ are
%attain
%second order accuracy in nominal coverage probabilities and
significantly more accurate than their counterparts based on normal
critical points. This observation is also corroborated
by the simulation results of Section~\ref{sec5}.

In Section~\ref{sec3-sub4}, a further refinement is obtained. A more
careful analysis of the $o_p(n^{-1/2})$-term in Theorem~\ref{thm-st}
shows that although it outperforms the normal approximation over
the class $\cC_q$, this rate does not always match the ``optimal''
level, namely $O_p(n^{-1})$ that is attained by the bootstrap
in the more classical
setting of estimation of regression parameters by the OLS
method with a \textit{fixed} $p$.
%The main reason behind this lack
%of 'optimality' of the bootstrap approximation for $\R_n$ is the
%effect of the bias, introduced by the penalization part
%in \eqref{alasso-def}.
Exploiting the higher order analysis in the proof of Theorem \ref
{thm-st}, we
carefully construct a modified studentized version $\breve{\R}_n$
of~$\bbln$. Theorem~\ref{thm3-3} shows that under
slightly stronger regularity
conditions (compared to those in Theorem~\ref{thm-st}),
the rate of bootstrap approximation
for the modified pivot $\breve{\R}_n$ is $O_p(n^{-1})$.
This appears to be a remarkable result because, even with a diverging
$p$ and with the regularization step,
the specially constructed pivotal quantity $\breve{\R}_n$
attains the same optimal
rate $O_p(n^{-1})$ as in the classical set up of
linear regression with a fixed $p$.
%An immediate practical implication of this is that for $q=1$,
% two-sided bootstrap confidence intervals (CIs) based on
% $\breve{\R}_n$ achieve the same level of accuracy as in the
% classical set up of bootstrapping the OLS of the regression
%parameters where $p$ is bounded and no penalization is
%used.

The key technical tool used in the proofs of the results in
Sections~\ref{sec3} and~\ref{sec4} is
an Edgeworth expansion (EE) result for the ALASSO estimator and
its studentized version, given in Theorem~\ref{thm5-1} of
Section~\ref{sec6}, which may be of independent interest.
The derivation of the EE
critically depends on the choice of
the initial estimator in \eqref{alasso-def}.
In Sections~\ref{sec3} and~\ref{sec4}, the initial estimator is
chosen to be the OLS, which necessarily
requires $p\leq n$. However, in many applications,
it is important to allow $p>n$. In such situations,
one may use a bridge estimator [cf. \citet{knight-fu}]
in place of the OLS as the initial estimator.
In Section~\ref{sec-new-th}, we show that under some suitable regularity
conditions, the bootstrap approximation to the distributions
of $\R_n$ and $\breve{\R}_n$ continue to be second order
correct even for $p>n$. Here, $p$ is allowed to grow at
polynomial rates in $n$. More precisely, we allow $p=O(n^a)$ for any given
$a>1$, provided (in addition to
certain other conditions) $\E|\ep_1|^r<\infty$
for a sufficiently large $r$, depending on $a$.
Thus, the allowable growth rate of $p$ depends on the
rate of decay of the tails of the error distribution.

The rest of the paper is organized as follows. We conclude this
section with a brief literature review. In Section~\ref{sec2}, we
introduce the theoretical framework and state the regularity
conditions. Results on the
rate of convergence to the oracle limit law is given in Section~\ref{sec3}.
%A brief description of the residual bootstrap method for
%the ALASSO and
The main results on the bootstrap
are given in Section~\ref{sec4} for the $p\leq n$ case
and in Section~\ref{sec-new-th} for the $p>n$ case.
Section~\ref{sec5} presents the results
from a moderately large simulation study
%illustrating the
%theoretical results from Sections~\ref{sec3}-\ref{sec-new-th},
and it also gives two real data examples. An outline of the proofs
of the main results is given in Section~\ref{sec6} and their
detailed proofs are relegated to a supplementary material file; cf.
\citet{CL}.

The literature on penalized regression in high dimensions
has been growing very rapidly in recent years; here we
give only a modest account of the work
that is most related to the present paper
due to space limitation. In two important papers, \citet{tibs-96}
%Tibshirani (1996)
introduced the LASSO,
as an estimation and
variable selection method
% in regression models with $\ell_1$-based
%penalization. In a seminal work,
and \citet{zou-06} introduced the ALASSO method
as an improvement over the LASSO and established its
oracle property.
% model-consistency and asymptotic distribution theory
%under weaker conditions on the design matrix.
Other popular penalized estimation and variable selection methods
are given by the SCAD [\citet{fan-li-01}]
and the Dantzig Selector [\citet{candes-tao-07}].
Properties of the ALASSO and the related methods
have been investigated by
many authors, including
\citet{knight-fu},
\citet{meinshausen-2006},
\citet{wainwright-2006},
\citet{bunea-t-w-07},
\citet{bickel-dantzig-2008},
\citet{huang-ma-sinica-08},
\citet{huang-ma-aos-08},
\citet{zhang-2008-36},
\citet{meinshausen-yu-09},
\citet{potscher-schneider-09},
\citet{ac-snl-sankhya-11},
\citet{gupta-12}
among others. \citet{fan-li-01} introduced the important
notion of ``oracle property'' in the context of penalized estimation
and variable selection by the SCAD. Post model selection inference,
including the bootstrap
and its variants have been investigated by
\citet{bach-09}, Chatterjee and
Lahiri (\citeyear{ac-snl-Lasso-08,ac-snl-jasa-11}),
\citet{minnier-jasa-11} and
\citet{berk-et-al-12},
among others.

%s2 #&#
%s2 ###
\section{Preliminaries and the regularity conditions}
\label{sec2}
%s2.1 #&#
%s2.1 ###
\subsection{Theoretical set up}
\label{sec2-sub1}
For deriving the theoretical results, we consider a generalized
version of \eqref{true-model}, where $p=p_n$ is allowed to
depend on the sample size $n$.
% Accordingly,
%the regression parameter and the design vectors $\xb_i$'s of
%dimension $p_n\times1$
%are also allowed to
%depend on $n$.
To highlight this, we shall denote the true
parameter value by $\bb_n$
and redefine
\[
\T_n = \sqrt{n}\D_n (\bbln- \bb_n ),
\]
where, as in Section~\ref{sec1}, $\D_n$ is a $q\times p_n$
(known) matrix satisfying $\operatorname{tr} (\D_n \D_n') =
O(1)$, and $q$ does not depend on $n$.
%Thus, $\T_n$
%is a vector of $q$ linear combinations of the centered and
%scaled ALASSO estimator $\sqrt{n}(\bbln-\bb_n)$.
Also, for the $p\leq n$ case, that is, in Sections~\ref{sec3} and~\ref{sec4}, we shall take the
initial estimator $\tilde{\bb}_n$ to be the OLS
of $\bb_n$, given by
$
\tilde{\bb}_n = [\sum_{i=1}^n \xb_i\xb_i^\prime]^{-1}
\sum_{i=1}^n \xb_iy_i$.

Let $I_n=\{j\dvtx1\leq j\leq p_n, \beta_{j,n}\neq0\}$ be the (population)
set of nonzero regression coefficients, where $\beta_{j,n}$ is
the $j$th component of $\bb_n$.
The ALASSO yields an estimator
$
\widehat{I}_n \equiv\{j\dvtx1\leq j\leq p_n,
\widehat{\beta}_{j,n}\neq0\}
$
of $I_n$. For notational simplicity, we shall assume that
$I_n=\{1,\ldots, p_{0n}\}$ and also suppress the dependence
on $n$ in $p_n$, $p_{0n}$, etc., when there is no
chance of confusion.

%For the sake of completeness, we now briefly describe
%the residual bootstrap which is the most common
%method for bootstrapping a linear regression model
%with a nonrandom $\X_n$ matrix (cf. Freedman (1981)).
%First, define the residuals $e_i = y_i - \xb^{\prime}_i \bbln$,
% $i=1,\ldots,n$,
%based on the ALASSO estimator. Let $\check{e}_i = e_i- \bar{e}_n$,
%$i=1,\ldots,n$ denote the centered residuals, where
%$\bar{e}_n = n^{-1}\sum_{i=1}^n e_i$. Next, select a
% random sample of size $n$ with
% replacement from the collection of all
%centered residuals $\{ \check{e}_1,\ldots,\check{e}_n\}$,
%and denote it by $\{e^{*}_1,\ldots,e^{*}_n\}$. Using these
%resampled residuals, generate the residual bootstrap observations
%$$
%y^{*}_i = \xb^{\prime}_i \bbln+ e^{*}_i, i = 1,\ldots,n.
%$$
%Note that the centering step ensures the model requirement $E
%for the bootstrap error variable $e^{*}_1$.
%%Based on the observations $\{(y^{*}_i, \xb^{\prime}_i)
%%:i=1,\ldots,n\}$, we define
%The bootstrap version of a statistic is defined by
%replacing $\{(y_i, \xb^{\prime}_i):i=1,\ldots,n\}$
%with $\{(y^{*}_i, \xb^{\prime}_i):i=1,\ldots,n\}$
%and $\bb_n$ with $\bbln$. For example,
%the bootstrap version ALASSO
% estimator is given by
%{\left(y^{*}_i-\xb^{\prime}_i\ub\right)}^2
%+\lambda_n\sum_{j=1}^p \frac{|u_j|}{{|\tilde{\beta}^{*}_{j,n}|}^
%where $\bbtnb= {\left(\tilde{\beta}^{*}_{1,n},\ldots,
%estimator of $\bb_n$, with $y^{*}_i$ replacing the $y_i$'s. Similarly,
%the bootstrap version of $\T_n$ is defined as
%$$
%$$
%%%%%%%%%%%%%%%%%%%%

%s2.2 #&#
%s2.2 ###
\subsection{Conditions}
\label{sec3-sub1}
Let $\C_n = n^{-1}\sum_{i=1}^n \xb_i\xb_i^\prime$.
Write
$\C_n= ( (c_{i,j,n} ) )$
and\break
$\C^{-1}_n= ( (c^{i,j}_n ) )$, when it exists.
Partition
$\C_n$ as
\[
\C_n=\left[ \matrix{
\C_{11,n}, \C_{12,n}
\vspace*{2pt}\cr
\C_{21,n}, \C_{22,n} }
\right],
\]
where $\C_{11,n}$ is $p_{0}\times p_{0}$. Similarly, let
$\D^{(1)}_n$ is the $q\times p_0$ submatrix of $\D_n$,
consisting of the first $p_0$
%all
columns of $\D_n$.
% indexed by the elements of $I_n$.
Let $\bar{\xb}_n = n^{-1}\sum_{i=1}^n \xb_i$ and let $\bar{\xb}^{(1)}_n$
denote the first $p_0$ components of $\bar{\xb}_n$. Define
\[
\bolds{\Sigma}^{(0)}_n  = %
\left[\matrix{
\D^{(1)}_n\C^{-1}_{11,n}{ \bigl(
\D^{(1)}_n \bigr)}{}^{\prime}\sigma^2 &
\D^{(1)}_n\C^{-1}_{11,n}\bar{
\xb}^{(1)}_n\cdot\E \bigl(\varepsilon^3_1
\bigr) \vspace*{2pt}
\cr
{ \bigl(\bar{\xb}^{(1)} \bigr)}{}^{\prime}
\C^{-1}_{11,n}{ \bigl(\D^{(1)}_n
\bigr)}^\prime\cdot\E \bigl(\varepsilon^3_1
\bigr) & \Var \bigl(\varepsilon^2_1 \bigr)} \right],
\]
which is used in
%the statement of
condition (C.3) below.
%plays an important role in the derivation of the EE
%for the studentized statistic $R_n$.
Let $\A_{i\cdot}$ and $\A_{\cdot j}$,
respectively, denote the $i$th row and the $j$th
column of a matrix $\A$, and let $\A^\prime$
denote the transpose of $\A$.
For $x, y\in\mathbb{R}$, let $x\vee y = \max\{x,y\}$,
$x_{+}=\max\{x,0\}$ and $\operatorname{sgn}(x) = -1,0,1$
according as $x<0$, $x=0$ and $x>0$.
Let $\iota= \sqrt{-1}$.
Unless otherwise stated, limits in
the order symbols are taken by letting $\nti$.

We shall make use of the following conditions:

\begin{enumerate}[(C.6)]
\item[(C.1)] There exists $\delta\in(0,1)$, such that
for all
$n>\delta^{-1}$,
\label{c1}
\[
{ \bigl(\xb^{\prime}\C_{12,n}\yb \bigr)}^2  \leq
\delta^2 \bigl(\xb^{\prime}\C_{11,n}\xb \bigr)\cdot
\bigl(\yb^{\prime}\C_{22,n}\yb \bigr)\qquad \mbox{for all $\xb\in
\mathbb{R}^{p_0}$, $\yb\in\mathbb{R}^{p-p_0}$.}
\]
\item[(C.2)] Let $\eta_n$ and $\eta_{11,n}$ denote
the smallest eigen-values of $\C_n$ and $\C_{11,n}$, respectively.
\label{c2}
\begin{enumerate}[(ii)]
\item[(i)] $\eta_{11,n}>Kn^{-a}$ for some $K\in(0,\infty)$ and
%{\footnote{doubt: what is the right endpoint for $a$?}}
$a\in[0,1]$.
\label{c2i}
\item[(ii)]
\label{c2ii}
$
\max\{n^{-1}\sum_{i=1}^n ({|x_{i,j}|}^r +
{|\tilde{x}_{i,j}|}^r)\dvtx1\leq j\leq p\} = O(1),
$
where $\tilde{x}_{i,j}$ is the $j$th element of
${ (\xb^{\prime}_i \C^{-1}_n )}$ (for $p\leq n$)
and $r\geq3$ is an integer
(to be specified in the statements of theorems).
\end{enumerate}
%
%%%%%%%%%%%%%%%%%%%%%old (C.6)
%
\item[(C.3)] There exists a $\delta\in(0,1)$ such that for all
$n>\delta^{-1}$:
\label{c3}
\begin{enumerate}[(ii)$^\prime$]
\item[(i)]
\label{c3i}
$
\sup\{\xb^{\prime}\D^{(1)}_n\C^{-1}_{11,n}{(\D^{(1)}_n)}^\prime
\xb\dvtx \xb\in\mathbb{R}^q, \|\xb\|=1\} < \delta^{-1}
$.
\item[(ii)]
\label{c3ii}
$
\inf\{\xb^{\prime}\D^{(1)}_n\C^{-1}_{11,n}{(\D^{(1)}_n)}{}^{\prime
}\xb\dvtx \xb\in\mathbb{R}^q, \|\xb\|=1\} > \delta.
$
\item[(ii)$^\prime$]
\label{c3iip}
$
\inf\{\tb^{\prime}\bolds{\Sigma}^{(0)}_n\tb\dvtx\tb\in\mathbb
{R}^{q+1}, \|
\tb\|=1\}>\delta
$.
\end{enumerate}
\item[(C.4)]
\label{c4}
$
\max\{|\beta_{j,n}|\dvtx j\in I_n\} = O(1)$
{and} $ \min\{|\beta_{j,n}|\dvtx j\in I_n\}\geq Kn^{-b},
$
for some $K\in(0,\infty)$ and $b\in[0,1/2)$, such that
$a+2b\leq1$, where $a$ is as in (C.2)(i):
%%%%%%%%%%%%%%%%%%%%Old (C.3)
%
\item[(C.5)]
\label{c5}
\begin{enumerate}[(ii)$^{\prime}$]
\item[(i)] $\E(\varepsilon_1)=0$, $\E(\varepsilon^2_1)
=\sigma^2\in(0,\infty)$ and $\E{|\varepsilon_1|}^r<\infty$, for some
$r\geq3$.
\label{c5i}
\item[(ii)] \label{c5ii}
$\varepsilon_1$ satisfies Cram\'{e}r's condition:
$
\limsup_{|t|\rai} |\E(\exp(\iota t\varepsilon_1) )| < 1.
$

\item[(ii)$^{\prime}$] \label{c5iip}
$(\varepsilon_1,\varepsilon^2_1)$ satisfies Cram\'{e}r's
condition,
\[
\limsup_{\|(t_1,t_2)\|\rai} \bigl|\E\exp \bigl(\iota\cdot \bigl(t_1
\varepsilon_1 + t_2 \varepsilon^2_1
\bigr) \bigr)\bigr | < 1.
\]
\end{enumerate}
%
%%%%%%%%%%%%%%%%old (C.5)
%
\item[(C.6)]
\label{c6}
%The regularization parameter $\lambda_n\in(0,\infty)$ and the
%tuning parameter $\gamma$ satisfy the following conditions:
There exists $\delta\in(0,1)$ such that for all $n\geq\delta^{-1}$,
\begin{eqnarray*}
\frac{\lambda_n}{\sqrt{n}} & \leq&\delta^{-1} n^{-\delta} \min \biggl\{
\frac{n^{-b\gamma}}{p_0},\frac{n^{-b\gamma-
{a}/{2}}}{\sqrt{p_0}}, n^{-a} \biggr\}\quad \mbox{and}
\\
\frac{\lambda_n}{\sqrt{n}}\cdot n^{\gamma/2} &\geq&\delta n^{\delta
} \max \bigl
\{n^{a}p_0, p^{3/2}_0
n^{b{(1-\gamma)}_{+}} \bigr\}.
\end{eqnarray*}
\end{enumerate}

We now comment on the conditions. Condition (C.1)
is equivalent to saying that the multiple correlation between
relevant variables (with $\beta_{j,n}\neq0$) and the spurious
variables ($\beta_{j,n}=0$) is strictly less than one, in
absolute value. This condition is weaker than assuming
orthogonality of the two sets of variables. Variants of this
condition has been used in the literature, particularly in
the context of the Lasso; see \citet{meinshausen-yu-09},
\citet{huang-ma-aos-08}, \citet{ac-snl-jasa-11}, and the
references therein.

Condition (C.2) gives the
regularity conditions on the design matrix that are needed for
establishing an $(r-2)$th order EE
%Edgeworth expansion
for the ALASSO estimator and its bootstrap versions.
(C.2)(i) requires a lower bound on the
smallest eigen-value of the submatrix $\C_{11,n}$ corresponding
to the relavent variables (with $\beta_{j,n}\neq0$), in
the increasing dimensional case. When $p$ is bounded,
$\C_n\rightarrow\C$ (elementwise) and~$\C$ is nonsingular,
this condition holds with $a=0$. Condition (C.2)(ii)
is a uniform bound on the $\ell_r$-norms of the sequences
$\{x_{i,j}\}^n_{i=1}$, $\{\tilde{x}_{i,j}\}^{n}_{i=1}$, that
are needed
for obtaining a uniform bound on the $r$th order moments of
the weighted sums $\sum_{i=1}^n x_{i,j}\varepsilon_i$ and
$\sum_{i=1}^n \tilde{x}_{i,j}\varepsilon_i$, for $1\leq j\leq p$.\vspace*{1pt}
Note that for $r=2$, the condition
$
\max\{n^{-1}\sum_{i=1}^n {|x_{i,j}|}^r\dvtx1\leq j\leq p\} = O(1)
$
is equivalent to requiring that the diagonal elements of the
$p\times p$ matrix $\C_n$ be uniformly bounded. Similarly, for $r=2$,
\begin{eqnarray*}
n^{-1}\sum_{i=1}^n {|
\tilde{x}_{i,j}|}^r & = &{ \bigl(\C^{-1}_n
\bigr)}_{j\cdot} \Biggl(n^{-1}\sum_{i=1}^n
\xb_i\xb^{\prime}_i \Biggr) { \bigl(
\C^{-1}_n \bigr)}_{\cdot j}
\\
& =& { \bigl(\C^{-1}_n \bigr)}_{j\cdot}
\C_n { \bigl(\C^{-1}_n \bigr)}_{\cdot j} =
{ (\mathbb{I}_p )}_{j\cdot}{ \bigl(\C^{-1}_n
\bigr)}_{\cdot j} = c^{j,j}_n,
\end{eqnarray*}
where $\mathbb{I}_p$ denotes the identity matrix of order $p$. Thus,
for $r=2$,
%
%e2.1 #&#
%
%e2.1 ###
\begin{eqnarray}
\max \Biggl\{n^{-1}\sum_{i=1}^n
{|\tilde{x}_{i,j}|}^r\dvtx1\leq j\leq p \Biggr\} &= O(1),
\label{eq3-1}
\end{eqnarray}
if and only if the diagonal elements of $\C^{-1}_n$ are uniformly bounded.
Condition~(C.2)(ii) is a stronger version of
these conditions with $r\geq3$, dictated by the order of the
EE
% Edgeworth expansion
one is interested in.
%
% Note that the condition
%holds if the $x_{i,j}$'s are bounded in the absolute value by
%some $K\in(0,\infty)$, for all $i,j$ and $n$.

%%%%%%%%%Discusssion of (C.3)
Conditions (C.3)(i) and (C.3)(ii)
require that the maximum and the minimum eigen-values of the
$q\times q$ matrix $\D^{(1)}_n\C^{-1}_{11,n}{(\D^{(1)}_n)}{}^\prime$
be bounded away from zero and infinity, respectively.
A sufficient condition is
the existence of a nonsingular
limit of $\D^{(1)}_n\C^{-1}_{11,n}{(\D^{(1)}_n)}{}^{\prime}$,
which we do \textit{not} assume. (C.3)(ii)$^\prime$
is a stronger form of (C.3)(ii) that is needed
for the studentized case only. Note that (C.3)
rules out inference on individual zero components
of $\bb_n$ (as $\D_n^{(1)}=\Zero$ in this case). The
main results of the paper are valid only for linear combinations
of the ALASSO estimator that put nontrivial weights on at least one nonzero
component of $\bb_n$.

Next consider condition (C.4)
which makes it possible to separate out
the signal from the noise by the ALASSO.
It requires the minimum of the nonzero coefficients to be of
coarser order than $O(n^{-1/2})$, so that the coefficients
are not masked by the estimation error, which is of the order
$O_p(n^{-1/2})$. It is worth pointing out that
the results of the paper remain valid
if the requirement $a+2b\leq1$ in condition (C.4)
is replaced by a somewhat weaker condition $n^{a+2b} = O(np_0)$.
Condition~(C.5)
is a moment and smoothness condition on the error variables. These
are required for the validity of an $(r-2)$th order EE, $r\geq3$, where
(C.5)(ii) is used for $\T_n$ and its stronger version
(C.5)(ii)$^\prime$ for the studentized cases, respectively.

%Edgeworth expansion.

Finally, consider condition (C.6). When $p_0$,
the number of nonzero components of $\bb_n$ is fixed
(but the total number of parameters $p$ may tend to $\infty$),
we may suppose that $\bb_n=\bb$ for all $n\geq1$ and hence, the
nonzero components of $\bb_n$ are bounded away from zero. If,
in addition, the submatrix $\C_{11,n}$ converges elementwise
to a $p_0\times p_0$ nonsingular matrix $\C$,
then $a=b=0$. In this case, condition~(C.6)
is equivalent to
\[
\frac{\lambda_n}{\sqrt{n}} + { \biggl[\frac{\lambda_n}{\sqrt{n}} \cdot n^{\gamma/2}
\biggr]}^{-1} = O \bigl(n^{-\delta} \bigr)
\]
for some $\delta>0$. This condition may be compared to the condition
\[
\frac{\lambda_n}{\sqrt{n}} + { \biggl[\frac{\lambda_n}{\sqrt{n}} \cdot n^{\gamma/2}
\biggr]}^{-1} = o(1),
\]
that was imposed by \citet{zou-06} to establish the asymptotic
distribution (and the {\textit{oracle}} property) of the ALASSO,
further assuming that $p$ itself is fixed. Thus, for a regression
problem with finitely many nonzero regression parameters and a
{\it{nice}} design matrix, the EE results hold under a
slight strengthening of the \citet{zou-06} conditions on
$\lambda_n$ and $\gamma$. It is interesting to note that the
growth rate of the zero components $(p-p_0)$ (or $p$ itself)
does not have a direct impact on $\lambda_n$ and $\gamma$ in
condition (C.6). However, when either $p_0\rai$\vadjust{\goodbreak}
or some of the nonzero components of $\bb_n$ become small,
the choices of $\lambda_n$ and $\gamma$ start to depend on the
associated rates. A similar behavior ensues for a nearly
singular submatrix $\C_{11,n}$. Further, note that for any
given values of $a\in[0,1]$ and $b\in[0,1/2)$,
we may allow $p_0=O(n)$ (with $p_0\leq n$), by choosing
$\lambda_n$ and $\ga^{-1}$ suitably small.
See Remark~\ref{rem-1}
in Section~\ref{sec3}
for more details on the
implications of
these conditions.\vspace*{-3.5pt}

%s3 #&#
%s3 ###
\section{Rates of convergence to the oracle distribution}
\label{sec3}
The main results of this section give upper and lower bounds on the
accuracy of approximation by the limiting oracle distribution for
the ALASSO. To describe the terms in the bounds,
let
% $\bar{\xb}^{(1)}_n = n^{-1}\sum_{i=1}^n \xb^{(1)}_i$
% and
$\bbf_n=\D^{(1)}_n\C^{-1}_{11,n}\sbf^{(1)}_n\cdot
\frac{\lambda_n}{\sqrt{n}},
$\vspace*{-1pt}
where $\sbf^{(1)}_n$ is a $p_0\times1$
vector with $j$th component
$
s_{j,n} = \operatorname{sgn}(\beta_{j,n}){{|\beta_{j,n}|}^{-\gamma}},
1\leq j\leq p_0
$.
Also let
$
\bolds{\Gamma}_n = \D^{(1)}_n\C^{-1}_{11,n}\bolds{\La
}^{(1)}_n\C
^{-1}_{11,n}{(\D^{(1)}_n)}^{-1}
$
where $\bolds{\La}^{(1)}_n$ is a diagonal matrix with $(j,j)$th element
given by $\operatorname{sgn}(\beta_{j,n}){|\beta_{j,n}|}^{-(\gamma
+1)}$, $1\leq j
\leq p_0$. Also, for a $k\times k$ nonnegative definite matrix
$\bolds
{\Sigma}$,
let $\Phi(\cdot\dvtx\bolds{\Sigma})$ denote the Gaussian measure on
$\mathbb{R}^k$
with zero mean and covariance matrix $\bolds{\Sigma}$.

Then we have the following result:\vspace*{-3.5pt}

%th3.1 #&#
%
\begin{theorem}
\label{thm3-1}
Suppose that conditions \textup{(C.1)--(C.6)}
hold with $r=4$ and that $\tilde{\bb}_n$
is the OLS of $\bb_n$. Then
\begin{eqnarray*}
\Delta_n & \equiv&\sup_{B \in\cC_q} \bigl|\pr(
\T_n\in B ) - \Phi \bigl(B\dvtx\sigma^2
\D^{(1)}_n\C^{-1}_{11,n}{ \bigl(
\D^{(1)}_n \bigr)}^{\prime} \bigr) \bigr|
\\[-3pt]
& =& O \biggl(n^{-1/2}+\|\bbf_n\|+\frac{\lambda_n}{{n}}\cdot
n^{a+b(\gamma+1)} \biggr).\vspace*{-6pt}
\end{eqnarray*}

\end{theorem}

Theorem~\ref{thm3-1} gives a precise description of the quantities
that determine the rate of convergence to the normal limit.
In particular,
%as noted by \citet{zou-06},
the ALASSO
estimator has a bias that may lead to an inferior rate
of convergence to the limiting normal distribution
[compared to the standard $O(n^{-1/2})$ rate],
depending on the choice of the penalty constant $\lambda_n$,
the exponent $\gamma$ and the rate\vspace*{1pt} of decay of the
smallest of the regression parameters. In addition, there
is a third term, of the order $a_{3,n} \equiv
\lambda_n \cdot n^{-1+a+b(\gamma+1)}$ that results from
the use of the initial estimator $\bbtn$ in the ALASSO
penalization scheme and that can also lead to
a sub-$n^{-1/2}$-rate of convergence to the normal limit.

We next show that under some mild conditions, the
bound given in Theorem~\ref{thm3-1} is precise in the
sense that, in general, it cannot be improved upon.\vspace*{-3.5pt}

%th3.2 #&#
%
\begin{theorem}
\label{thm3-1b}
Suppose that the conditions of~\ref{thm3-1}
hold and that $\E\varepsilon_1^3 \neq0$,
$\liminf_{\nti} \sum_{|\alb|=3} |{ (\D^{(1)}_n
\C^{-1}_{11,n}\bar{\xb}^{(1)}_n )}^{\alb} |\neq0$,
$n^{a+b(\gamma+1)}=$
$O (\operatorname{tr}(\bolds{\Gamma}_n) )$
and $n^{b\gamma}=O (\|\D^{(1)}_n\C^{-1}_{11,n}\sbf^{(1)}_n\|
)$. Then
\[
\Delta_n  \asymp \biggl[n^{-1/2}+\frac{\lambda_n}{\sqrt{n}}\cdot
n^{b\gamma} + \frac{\lambda_n}{{n}}\cdot n^{a+b(\gamma+1)} \biggr],\vspace*{-2pt}
\]
where we write $a_n\asymp b_n$ if $a_n=O(b_n)$ and $b_n = O(a_n)$ as
$\nti$.\vadjust{\goodbreak}
\end{theorem}

Note that under the additional conditions of Theorem~\ref{thm3-1b},
the co-efficients of the first and the third terms on the right-hand\vspace*{1pt}
side of the display above
are
nonnegligible in the limit and $\|\bbf_n\|\geq K
\frac{\lambda_n}{\sqrt{n}}\cdot n^{b\gamma}$ for some
constant $K\in(0,\infty)$. As a result, the leading terms in the EE for
$\T_n$ that determine the upper bound in Theorem~\ref{thm3-1}
are also bounded from below by constant multiples of the three
factors appearing in Theorem~\ref{thm3-1b}. As a consequence,
the exact rate of approximation by the oracle distribution
to the centered and scaled ALASSO estimator $\T_n$
is given by the maximum of these three terms. In Remark~\ref{rem-1} below,
we discuss in more details the effects of
the choices of the penalty constant $\lambda_n$,
the exponent $\gamma$, etc.
on the accuracy of the oracle based normal approximation.

%
%re1 #&#
%
\begin{remark}
\label{rem-1}
Suppose that $\lambda_n\sim Kn^{c}$ for some $K\in(0,\infty)$
and
%{\footnote{doubt: should it be that $c>0$ ?}}
$c\in\mathbb{R}$ and let $\|\C_{11,n}^{-1/2}\sbf^{(1)}_n\| =
O(n^{\gamma b})$. Then $\|\bbf_n\| \leq\|\D_n^{(1)} \C
_{11,n}^{-1/2}\|
\cdot\|\C_{11,n}^{-1/2}\sbf^{(1)}_n\| \lambda_n/\sqrt{n}
= O(\lambda_n n^{-{1}/{2} + \gamma b})$. Hence,
under the conditions of Theorem~\ref{thm3-1},
the rate of normal approximation for $\T_n$
is given by
\[
\max \bigl\{n^{-1/2}, n^{c+b\gamma- 1/2}, n^{a+b(\gamma
+1)+c-1} \bigr\}.
\]
Here, a sub-optimal rate results if either
$b\gamma+ c>0$ or $a+b(1+\gamma)+c>1/2$. Further,
the bias term is the leading sub-optimal term
whenever
%
%e3.1 #&#
%
%e3.1 ###
\begin{equation}
%c+b\gamma-1/2>a+b(\gamma+1)+c-1,
a+b<1/2 \quad\mbox{and}\quad b\gamma+ c>0. \label{rem-eq1}
\end{equation}
In this case, using the EE results from Section~\ref{sec6}
[cf. Theorem~\ref{thm5-1}(a)], one can conclude that,
for a linear function of $\bb_n$ (i.e., for a $1\times p$ vector $\D_n$
with $q=1$), the errors in coverage probabilities of
\textit{both one and two-sided} confidence intervals (CIs)
based on the oracle normal critical points are
$O (n^{-{1}/{2}+(b\gamma+ c)} )$. This rate is much worse
than the available optimal rates, particularly in the two-sided case.

By a similar reasoning, the third term is the dominant
sub-optimal term whenever
%
%e3.2 #&#
%
%e3.2 ###
\begin{equation}
a+b>1/2 \quad\mbox{and}\quad a+b(\gamma+1)+c\in(1/2,1). \label{rem-eq2}
\end{equation}
In this case, Theorem~\ref{thm5-1}(a) shows that one-sided
CIs based on the {oracle} distribution
r has a sub-optimal error. However, as the corresponding term in
the EE for $\T_n$ is even, it no longer contributes to
the error of coverage probability in the two-sided case.

Finally the optimal rate of convergence in Theorem~\ref{thm3-1b} holds,
provided
\[
c+b\gamma\leq0 \quad\mbox{and}\quad a+b(\gamma+1)+c\leq1/2.
\]
Since $a\geq0$, $b\geq0$ and $\gamma>0$, the first inequality
requires $c\leq0$, that is, $\lambda_n=O(1)$. Further,
for $ab>0$, that is, when
both the smallest eigen-value $\eta_{11,n}$ of $\C_{11,n}$
and the minimum of the nonzero components
(say $\beta_{1n}^{\mathrm{min}}$) of the regression vector
$\bb_n$ tend to {zero}, these inequalities require that
$c $ be chosen to be a sufficiently big
negative number (and thus, $\lambda_n$ to be a \textit{small} positive
number). This in turn leads to an inferior performance of the
ALASSO for variable selection. In the next section, we show that
the bootstrap attains the optimal rate of approximation to the
distribution of $\T_n$ without requiring such unreasonable conditions
on the choice of~$\lambda_n$.
\end{remark}

%s4 #&#
%s4 ###
\section{Accuracy of the bootstrap}
\label{sec4}

%s4.1 #&#
%s4.1 ###
\subsection{The residual bootstrap}
\label{sec2-sub2}
For the sake of completeness, we now briefly describe
the residual bootstrap
%which is the most common
%method for bootstrapping a linear regression model
%with a nonrandom $\X_n$ matrix
[cf. \citet{freedman-81}].
Let $e_i = y_i - \xb^{\prime}_i \bbln$, $i=1,\ldots,n$
denote the residuals
based on the ALASSO estimator, and let $\check{e}_i = e_i- \bar{e}_n$,
$i=1,\ldots,n$, where
$\bar{e}_n = n^{-1}\sum_{i=1}^n e_i$. Next, select a
random sample of size $n$ with
replacement from $\{ \check{e}_1,\ldots,\check{e}_n\}$,
and
denote it by $\{e^{*}_1,\ldots,e^{*}_n\}$.
Define the residual bootstrap observations
\[
y^{*}_i = \xb^{\prime}_i \bbln+
e^{*}_i, \qquad i = 1,\ldots,n.
\]
Note that the centering step ensures the model requirement $\E
\varepsilon
_1 =0$
for the bootstrap error variable $e^{*}_1$.
%Based on the observations $\{(y^{*}_i, \xb^{\prime}_i)
%:i=1,\ldots,n\}$, we define
The bootstrap version of a statistic is defined by
replacing $\{(y_i, \xb^{\prime}_i)\dvtx i=1,\ldots,n\}$
with $\{(y^{*}_i, \xb^{\prime}_i)\dvtx i=1,\ldots,n\}$
and $\bb_n$ with~$\bbln$. For example,
the bootstrap version ALASSO
estimator is given by
%
%e4.1 #&#
%
%e4.1 ###
\begin{equation}
\bolds{\beta}^{*}_n= \mathop{\operatorname{argmin}}_{\ub\in
\mathbb{R}^p}
\sum_{i=1}^n { \bigl(y^{*}_i-
\xb^{\prime}_i\ub \bigr)}^2 +\lambda_n
\sum_{j=1}^p \frac{|u_j|}{{|\tilde{\beta
}^{*}_{j,n}|}^\gamma},
\label{boot-alasso-def}
\end{equation}
where $\widetilde{\bolds{\beta}}^{*}_n= { (\tilde{\beta
}^{*}_{1,n},\ldots,
\tilde{\beta}^{*}_{p,n} )}{}^{\prime}$ is the bootstrap version
of the
initial estimator~$\tilde{\bb}_n$
(which is given by the OLS in this section),
obtained by replacing the $y_i$'s with $y^{*}_i$'s.
The bootstrap version of $\T_n$ is then defined as
$
\T^{*}_n = \sqrt{n}\D_n (\bolds{\beta}^{*}_n- \bbln)
$. Similarly, define $\R_n^*$ and $\breve{\R}_n^*$.

%s4.2 #&#
%s4.2 ###
\subsection{Rates of bootstrap approximation for $\T_n$}
\label{sec4-sub1}
The following result shows that the bootstrap approximation to
the distribution of $\T_n$ attains the rate $O_p (n^{-1/2} )$
under regularity conditions (C.1)--(C.6).

%th4.1 #&#
%
\begin{theorem}
\label{thm3-2}
If conditions \textup{(C.1)--(C.6)} hold with
$r=4$, then
\[
\sup_{B\in\cC_q} \bigl|\ps \bigl(\T^{*}_n\in B
\bigr) - \pr(\T_n\in B )\bigr | = O_p \bigl(n^{-1/2}
\bigr).
\]
\end{theorem}

A comparison of Theorem~\ref{thm3-2} and the results of Section~\ref{sec3}
shows that the bootstrap approximation attains the optimal
rate $O_p (n^{-1/2} )$, \textit{irrespective} of the order of magnitudes
of the bias term $\|\bbf_n\|$ and of the third term $a_{3,n}$
in Theorem~\ref{thm3-1}. In particular, this rate is attainable
even when the smallest eigen-value $\eta_{11,n}$ of $\C_{11,n}$
or the minimum of the nonzero components
(say $\beta_{1n}^{\mathrm{min}}$) of the regression vector
$\bb_n$ tend to {zero}.
% solely under
% conditions \hyperref[c1]{(C.1)} - \hyperref[c6]{(C.6)}.
Most importantly, the bootstrap approximation to the
ALASSO estimator attains the same level of accuracy
in increasing dimensions as in the simpler
case of the OLS of regression parameters when the dimension
$p$ of the regression parameter is fixed and no
penalization is used. Thus, the bootstrap approximation
for $\T_n$ is in a way immune to
the effects of high dimensions.
% handles the complexities arising from high dimensionality of the
%problem with equal facility.

%s4.3 #&#
%s4.3 ###
\subsection{Rates of bootstrap approximation for $\R_n$}
\label{sec4-sub3}

As is well known in the fixed~$p$ case [cf. \citet{hall-book-92}],
the bootstrap gives a more accurate approximation when it
is applied to a pivotal quantity, such as a studentized version
of a statistic, rather than to its nonpivotal
version, like $\T_n$. Here we consider
the following studentized version of the ALASSO estimator:
\[
\R_n = \T_n/\widehat{\sigma}_n,
\]
where $\widehat{\sigma}_n^2 = n^{-1}\sum_{i=1}^n \check{e}_i^2$
and $\check{e}_1, \ldots,\check{e}_n$ are
the centered residuals (cf. Section~\ref{sec2-sub2}).
As explained in Section~\ref{sec1}, this differs from the
standard version of the studentized statistic
$
\tilde{\R}_n = \widehat{\V}_n^{-1/2} \T_n
$
where $\widehat{\V}_n$ is an estimator of the
asymptotic covariance matrix
$\V_n =
\sigma^2
\D^{(1)}_n\C^{-1}_{11,n}{(\D^{(1)}_n)}{}^{\prime}$
of $\T_n$ given by the
oracle limit distribution; cf. Theorem~\ref{thm3-1}.
Note that this studentized version of $\T_n$ can
be computationally highly demanding, particularly
for repeated bootstrap computation,
when $p_0$ is large. In comparison, the proposed studentized
version of $\T_n$ that we consider here is based only on
a scalar factor and hence, computationally
simpler.

The following result gives the rate of bootstrap approximation
to the distribution of $\R_n$. For notational compactness,
in the rest of this section,
we shall write
(C.1)$^\prime$--(C.6)$^\prime$,
%
%(C.1)$^\prime$-(C.6)$^\prime$
to denote conditions (C.1)--(C.6), when
(C.3) and (C.6)
are defined with part~(ii)$^\prime$ instead of part (ii).

%th4.2 #&#
%
\begin{theorem}
\label{thm-st}
If conditions \textup{(C.1)}$^\prime$--\textup{(C.6)}$^\prime$
hold with $r=6$, then
\[
\sup_{B\in\cC_q} \bigl|\ps \bigl(\R^{*}_n\in B
\bigr) - \pr(\R_n\in B ) \bigr|  = o_p \bigl(n^{-1/2}
\bigr).
\]
\end{theorem}

Theorem~\ref{thm-st} shows that under conditions
(C.1)$^\prime$--(C.6)$^\prime$,
the bootstrap approximation to the distribution of $\R_n$
is second-order-correct, as it corrects for the
effects of the leading terms in the EE of $\R_n$.
From the proof of Theorem~\ref{thm5-1}, it follows
that the bootstrap not only captures the
usual $O(n^{-1/2})$ term in the EE, but it also
corrects for the effects of the second
and the third terms in the upper bound
of Theorem~\ref{thm3-1} that result from the
penalization step in the definition of the
ALASSO. The accuracy level $o_p(n^{-1/2})$ for the bootstrap
holds even when the actual magnitudes of these terms are
coarser than\vadjust{\goodbreak} $n^{-1/2}$ which, in turn, leads to a
poor rate of approximation by the limiting
normal distribution. A practical implication of this result
is that percentile-$t$ bootstrap CIs based on $\R_n$
will be more accurate than the CIs based on the large sample
normal critical points. Indeed, the finite sample
simulation results presented in Section~\ref{sec5} show
that the CIs based on normal critical points
are practically useless in moderate samples
and improvements in the coverage accuracy
achieved by the bootstrap CIs based on $\R_n$
are %strikingly
spectacular.\vspace*{-3pt}

%s4.4 #&#
%s4.4 ###
\subsection{A modified pivot and higher order correctness}
\label{sec3-sub4}
Although the residual bootstrap approximation for the studentized
statistic $\R_n$ is second order correct, a more careful analysis
shows that it may fail to achieve the
same {\it{optimal}} rate, namely,
$O_p(n^{-1})$ as in the traditional fixed
and finite dimensional regression problems.
The main reason behind this is the effect
of the bias term
% for symmetric sets in $\cC_q$ (cf. \citet{hall-book-92}),
$\|\bbf_n\|$ in Theorem~\ref{thm3-1}, which can be coarser
than $n^{-1/2}$. While the second order correctness is a desirable
property for the one-sided CIs, the higher level
of accuracy, namely $O_p(n^{-1})$, is important for
two-sided CIs;
cf. \citet{hall-book-92}.
%
% as the corresponding bootstrap confidence regions will
%have 'sub-optimal' coverage accuracy.
To that end,
we now define a modified pivotal quantity
%
%e4.2 #&#
%
%e4.2 ###
\begin{equation}
\breve{\R}_n  = \frac{\sqrt{n}\D_n (\bbln-
\bb_n ) + \breve{\bbf}_n}{\breve{\sigma}_n}, \label{rndef-alt}
\end{equation}
where $\breve{\bbf}_n = \breve{\D}^{(1)}_n{\breve{\C}^{-1}_{11,n}}
\breve{\sbf}^{(1)}_n\cdot\frac{\lambda_n}{\sqrt{n}}$, $\breve{\D
}^{(1)}_n$
and $\breve{\C}^{(1)}_{11,n}$ are, respectively, $q\times|\widehat{I}_n|$
and $|\widehat{I}_n|\times|\widehat{I}_n|$ submatrices of $\D_n$
and $\C_n$ with columns (and also rows, in case of $\breve{\C}_{11,n}$)
in $\widehat{I}_n=\{j\dvtx1\leq j\leq p, \widehat{\beta}_{j,n}\neq
0\}$,
and similarly, $\breve{\sbf}^{(1)}_n$ is the $|\widehat{I}_n|\times1$
vector with $j$th element $\operatorname{sgn}
(\widehat{\beta}_{j,n}){|\tilde{\beta}_{j,n}|}^{-\gamma}$,
$j\in\widehat{I}_n$. Here $\breve{\sigma}^2_n$ is
defined as
\[
\breve{\sigma}^2_n = \frac{1}{n}\sum
_{i=1}^n { (\breve{\varepsilon}_i -
\bar{\breve{\varepsilon}}_n )}^2,
\]
where
$
\breve{\varepsilon}_i = y_i - \xb^{\prime}_i\breve{\bb}_n,
\mbox{and}
\breve{\beta}_{j,n} = \tilde{\beta}_{j,n}\cdot
\mathbf{1}(j\in\widehat{I}_n), 1\leq j\leq p
$.
Note that $\breve{\R}_n$
%we use components of the Adaptive Lasso estimator in the sign of
% $\beta_{j,n}$'s and those of the preliminary estimator
%$\widetilde{\bb}_n$ to estimate the denominators.
is obtained by applying a specially designed
bias-correction term to $\T_n$ and by a
suitable rescaling, which are suggested by the
form of the third order EE of Theorem~\ref{thm5-1}. Also,
it is interesting to note that for both of these
estimators, we only use the sub-vectors
of the design vectors $\bfx_i$'s and components of the
initial estimator that correspond to the (random) set of
variables selected by the ALASSO. Next, define $\breve{\R}^{*}_n$,
the bootstrap version of $\breve{\R}_n$, by replacing $\{y_1,\ldots,y_n\}$
and $\bb$ by $\{y^{*}_1,\ldots,y^{*}_n\}$ and $\bbln$, respectively.
Then we have the following result:\vspace*{-3pt}

%th4.3 #&#
%
\begin{theorem}
\label{thm3-3}
If conditions \textup{(C.1)}$^{\prime}$--\textup{(C.6)}$^{\prime}$ hold
with $r=8$, then
\[
\sup_{B\in\cC_q} \bigl| \ps \bigl(\breve{\R}^{*}_n
\in B \bigr) - \pr(\breve{\R}_n\in B ) \bigr|  = O_p
\bigl(n^{-1} \bigr).
\]
\end{theorem}

Theorem~\ref{thm3-3} asserts that under appropriate regularity
conditions, the rate of bootstrap approximation to\vadjust{\goodbreak}
the modified pivotal quantity $\breve{\R}_n$ attains the
the ``optimal'' level of accuracy irrespective of
the magnitude of $\|\bbf_n\|$. An immediate consequence
of this result is that symmetric bootstrap confidence regions
based on the modified pivot attains the higher
rate $O(n^{-1})$ of convergence accuracy even when
the magnitude of $\|\bbf_n\|$ is coarser than $n^{-1/2}$.
As explained in Remark~\ref{rem-1}, the coarser
magnitude of $\|\bbf_n\|$ can occur quite naturally
in a variety of situations whenever
a combination of values of the underlying regression parameters,
the design matrix and the choice of the penalty constant satisfy
\eqref{rem-eq1}.
In such cases, bootstrap CIs based on $\breve{\R}_n$
gives a marked improvement
over normal critical points based
CIs where the accuracy is sub-$O(n^{-1/2})$
for both
one- and two-sided CIs.

%s5 #&#
%s5 ###
\section{Results for the $p>n$ case}
\label{sec-new-th}

In many applications, $p$ is much larger than~$n$, and post variable selection
inference on the regression parameters is an even more challenging
problem. In this section, we study properties of the bootstrap
approximation to the studentized ALASSO estimator in the
$p>n$ case.
Note that for $p>n$,
the $p\times p$ matrix $n^{-1}\sum_{i=1}^n \xb_i\xb^{\prime}_i$ is
always singular and hence the OLS of $\bb_n$ is no longer uniquely
defined. In the literature, a popular choice of
the
initial root-$n$ consistent estimator $\bbtn$
for $p>n$ is
the LASSO estimator, although other bridge estimators of $\bb_n$
[cf. \citet{knight-fu}]
can also be used.
Let $\bbln$ be the ALASSO estimator defined by \eqref{alasso-def},
with a root-$n$ consistent initial estimator $\bbtn$. Also
define the studentized version of
$\hat{\bb}_n$ (cf. Section~\ref{sec4-sub3}) by
$
\R_n = \widehat{\sigma}^{-1}_n \T_n$
where $\widehat{\sigma}^2_n$ is the average of squared centered
residuals $\breve{e}_1,\ldots,\breve{e}_n$, from the ALASSO fit,
and
define the bias corrected version $\breve{\R}_n$ as
in \eqref{rndef-alt}.

To prove the results in the $p>n$ case, we need
the following condition:

(C.7) There exists $K\in(0,\infty)$ such that
%
%e5.1 #&#
%
\begin{eqnarray}
\label{hd-cond} \pr \Bigl(\max_{1\leq j\leq p} \bigl| \sqrt{n} (\widetilde{
\beta}_{j,n} - \beta_{j,n} ) \bigr|>K \sqrt{\log{n}} \Bigr)& = &o
\bigl(n^{-1/2} \bigr),
\nonumber
\\[-8pt]
\\[-8pt]
\nonumber
\ps \Bigl(\max_{1\leq j\leq p} \bigl| \sqrt{n} \bigl(\widetilde{
\beta}^{*}_{j,n} - \widehat{\beta}_{j,n} \bigr) \bigr|>K
\sqrt{\log{n}} \Bigr) &=& o_p \bigl(n^{-1/2} \bigr).
\nonumber
\end{eqnarray}
We also need the following modified version of (C.2)(ii):

(C.2)(ii)$^\prime$
\[
\max_{1\leq j\leq p} \Biggl\{n^{-1} \sum
_{i=1}^n {|x_{i,j}|}^r \Biggr\}
+ \max_{1\leq j\leq{p_0}} \bigl\{c^{j,j}_{11,n} \bigr\} =
O(1),
\]
where $c^{j,j}_{11,n}$ is the $(j,j)$th element of $\C^{-1}_{11,n}$.

We now briefly discuss the conditions. Condition (C.7) is a high-level
condition that requires
the initial estimator $\bbtn$
and its bootstrap version not only to be $\sqrt{n}$-consistent,
but also to satisfy a suitable form of moderate deviation bound. For
estimators $\bbtn$, such that $\sqrt{n} (\widetilde{\beta
}_{j,n}-\beta_{j,n} )$ can be closely approximated by\vadjust{\goodbreak} $\sum_{i=1}^n
h_{j,i,n}\varepsilon_i$ for some $\{h_{j,i,n}\}\subset\mathbb{R}$ with
$\sum_{i=1}^n h^2_{j,i,n} = O(1)$, {(C.7)} holds if $\E{\varepsilon
^4_1}<\infty$ and $\sum_{i=1}^n h^4_{j,i,n}=o (n^{-1/2} )$.
See Proposition 8.4 [\citet{CL}] for an example. Condition
(C.2)(ii)$^\prime$ drops the condition
$
\max\{n^{-1} \sum_{i=1}^n {|\tilde{x}_{i,j}|}^r\dvtx1\leq j\leq p\}
= O(1),
$
in (C.2)(ii), which can no longer hold in the $p>n$ case,
as $\C^{-1}_n$ does not exist. Instead, it requires existence
of $\C^{-1}_{11,n}$, which is of dimension $p_0\times p_0$.
Thus, we must have $p_0\leq n$ (in addition to other conditions)
for the validity of the results in the $p>n$ case.

Let ${\R}^{*}_n$ and $\breve{\R}^{*}_n$
denote the (residual) bootstrap versions of $\R_n$ and
$\breve{\R}_n$, respectively. Then, we have the following result:

%th5.1 #&#
%
\begin{theorem}\label{thmhd-1}
Suppose that $p>n$ and conditions \textup{(C.1)},
\textup{(C.2)(i)},\break
\textup{(C.2)(ii)}$^\prime$, \textup{(C.3)--(C.7)} hold with $b=0$. Then
\begin{eqnarray*}
\sup_{B\in\mathcal{C}_q}\bigl|\pr(\R_n \in B) - \ps \bigl(
\R^{*}_n \in B \bigr)\bigr| & =& o_p
\bigl(n^{-1/2} \bigr)\quad \mbox{and}
\\
\sup_{B\in\mathcal{C}_q}\bigl|\pr(\breve{\R}_n \in B) - \ps \bigl(
\breve{\R}^{*}_n \in B \bigr)\bigr| & = &o_p
\bigl(n^{-1/2} \bigr).
\end{eqnarray*}
\end{theorem}

Thus, under the conditions of Theorem~\ref{thmhd-1}, the bootstrap
approximations based on the pivots $\R_n$ and $\breve{\R}_n$ are
both second-order accurate, even in the case where $p>n$. In comparison,
the oracle based normal approximation admits the sub-optimal bounds
of Section~\ref{sec3}, and therefore, it is significantly less accurate
than the bootstrap
approximations. This conclusion is also supported
by the finite sample simulation results
of Section~\ref{sec5} for the $p>n$ cases considered therein.

%A second restrictive condition in Theorem~\ref{thmhd-1} is $b=0$,
% which requires the nonzero coefficients to be bounded away from
% zero. Finally,

%re2 #&#
%
\begin{remark}
\label{rem-2}
Note that
in Theorem~\ref{thmhd-1},
the bound on the accuracy of the bootstrap
approximations to $\breve{\R}_n$ is just $o_p(n^{-1/2}$)
for the $p>n$ case. This is not as precise as
the bound in the $p\leq n$ case where it is $O_p(n^{-1}$).
It would be possible to derive a similar
bound for the $p>n$ case for $\breve{\R}_n$ if
we are willing to make some strong
additional assumptions on the initial estimator [e.g.,
existence of an EE for the joint
distribution of $\T_n$, $n^{-1}\sum_{i=1}^n {(\varepsilon^k_i -
\E\varepsilon_i^k)}$, with $k=1,2$ and suitable linear
combinations of $\sqrt{n}(\bbtn-\bb_n)$, which are
not known at this stage]. As a result, we do not pursue
such refinements here.\looseness=-1
\end{remark}

%re3 #&#
%
\begin{remark}
\label{rem-3}
Although we do not explicitly impose any growth conditions on~$p$ as a
function of $n$, there is, however, an implicit requirement
through condition (C.7). Indeed, if the leading terms in
$\sqrt{n}(\tilde{\beta}_{j,n}-\beta_{j,n})$ can be expressed
as $\sum_{i=1}^n h_{ji,n}\varepsilon_i$ for some $h_{j1,n},\ldots,h_{jn,n}\in
\mathbb{R}$ with $\sum_{i=1}^n h^2_{ji,n} = O(1)$, then
for~(C.7) to hold, arguments
in the proof of Lemma~\ref{lem5-1}(iii) require
that, for some integer $r\geq3$,
$\E{|\varepsilon_1|}^r<\infty$ and $p \cdot n^{-(r-2)/2}
= o(n^{-1/2})$. This implies that
%, $p = o(n^{(r-3)/2})$ when $\E{|\varepsilon_1|}^r<\infty$,
%or in other words,
$p$
can grow at a polynomial rate $p\sim Kn^a$, for some $K>0$ and $a>1$,
provided $\E{|\varepsilon_1|}^r<\infty$ for some $r > 2a+3$.
Thus, the allowable growth rate of $p$ depends on the lightness of the
tails of the error distribution.\vadjust{\goodbreak}
\end{remark}

%re4 #&#
%
\begin{remark}
\label{rem-4}
As pointed out by a referee, the use of $\bbtn$ in place of
$\widetilde{\bolds{\beta}}^{*}_n$ in
the bootstrap computation of the ALASSO estimator in \eqref
{boot-alasso-def} will yield a computationally more efficient
algorithm. It can be shown that with this modification, conclusions of
Theorems~\ref{thm-st},~\ref{thm3-3} and~\ref{thmhd-1} remain valid,
with the error bound $o_p(n^{-1/2})$ only.
\end{remark}

%s6 #&#
%s6 ###
\section{Simulation results}
\label{sec5}

%It is well known (cf. \citet{hall-book-92}) that under the 'smooth
%function model', the two-sided CIs for $1$-dimensional %parameters
%based on normal critical points have coverage accuracy $O(n^{-1})$,
%and this rate remains the same for the %bootstrap-$t$ CIs based on the
%Studentized statistic. In such cases, one uses the symmetric CIs (cf.
%i.e.} to the order $O(n^{-2})$. Here we show that even in the
%penalized %regression problem with a diverging number of parameters,
%symmetric bootstrap CIs based on the modified pivot $\breve{R}_n$
%%yield the same level of accuracy. Specifically, let $q=1$ and let $
%$\D_n$ is $1\times p$ is now a row-vector). For a given $\alpha
%the bootstrap distribution of $|\breve{R}_n|$. Set $\widehat{\theta}_n
%= \D_n\bbln$. Define the
%$2$-side d 100 $(1
% -\alpha)\%$ symmetric bootstrap CI for $\theta_n$, based on $
%$$
%J_n(\alpha) = \left[\widehat{\theta}_n - \frac{\breve{\bbf}_n+
%$$
%
%%%%%%%%%%%%5{\sc More to be added... }

In this section we study the finite sample performance of the proposed
bootstrap methods. The following cases corresponding to different
choices of $\bb_n$ were studied:
\begin{longlist}[(a)]
\item[(a)] $(n,p)=(60,10)$: with $p_0=5$ and $\bb_n = {(4,-1.5,-8,
0.9,-3, 0, \ldots, 0)}^\prime$.
\item[(b)] $(n,p) = (60,100)$: with $p_0 = 5$ and $\bb_n$ same as in
case (a) above, except that last 95 components are zeros.
\item[(c)] $(n,p) = (200, 80)$: with $p_0 = 10$ and
with the last 70 components being zeros,
\[
\bb_n = {(4, 2.5, 0.8, -1.5, -2, -5, -7.5, 5, 1.5, -3, 0,\ldots,0)}^\prime.
\]
\item[(d)] $(n,p)=(200,500)$: with $p_0 = 10$ and $\bb_n$ same as in
case (c) above, except that the last 490 components are zeros.
\end{longlist}
Cases (b) and (d) correspond to the $p>n$ case. In all cases, the
design vectors ${(x_{i,1},\ldots,x_{i,p_0})}^\prime$ are independently
generated from a normal population with mean~$\Zero$ and covariance
matrix $((\eta_{i,j}))$ with $\eta_{i,j}={(0.3)}^{|i-j|}$ and the
remaining $(p-p_0)$ covariates are i.i.d. $N(0,1)$. The errors $\{
\varepsilon_i\}$ are i.i.d. $N(0,1)$. We fix $\gamma= 1$. In the
high-dimensional case, since there is no unique least squares
estimator, we have used the LASSO estimator as the initial estimator
$\bbtn$, with associated tuning parameter $\lambda_{1,n}$. In the
ALASSO step, the penalty parameter is $\lambda_{2,n}$ and to avoid
division by zero, we used weights ${(|\widetilde{\beta
}_{j,n}|+a_n)}^{-1}$ with $a_n = n^{-1/2}$, to define the weighted
$\ell
_1$ penalty in \eqref{alasso-def}.

%s6.1 #&#
%s6.1 ###
\subsection{Comparison of oracle based normal CIs and bootstrap CIs}
\label{sim-sub0}

As suggested from Table~\ref{tab1}, in all cases
when the underlying true parameter value is large enough,
the bootstrap based CIs clearly superior to the oracle
based method. For moderately small underlying true parameters,
results in Table~\ref{tab2} suggest that the bootstrap-based methods
are still better than the Oracle method
for both one and two-sided CIs, even when $p>n$.
The improvement is most significant for the 2-sided CIs.

%t1 ###
\begin{table}
\caption{Comparison of empirical coverage probabilities and
average lengths (in parentheses) for~$90\%$ CIs for the underlying
parameter $\beta_1(=4)$ in cases \textup{(a)--(d)}. In all cases $\lambda
_{2,n} = 2n^{1/4}$ and in cases \textup{(b)} and \textup{(d)}, $\lambda_{1,n} = 0.5
n^{1/2}$}\label{tab1}
\begin{tabular*}{\textwidth}{@{\extracolsep{4in minus 4in}}lcccccc@{}}
\hline
& \multicolumn{3}{c}{\textbf{One-sided}} & \multicolumn{3}{c@{}}{\textbf{Two-sided (with average
lengths)}}\\[-4pt]
& \multicolumn{3}{c}{\hrulefill} & \multicolumn{3}{c@{}}{\hrulefill}\\
\textbf{Case} & $\bolds{\R_n}$ & $\bolds{\breve{\R}_n}$ & \textbf{Oracle} & $\bolds{\R_n}$ & $\bolds{\breve{\R}_n}$ &
\multicolumn{1}{c@{}}{\textbf{Oracle}}\\ %& $|\R_n|$\tmark[a] & $|\breve{\R}_n|$\tmark[a]
\hline
(a) & 0.898 & 0.904 & 0.668 & 0.918 & 0.900 & 0.158 %& 0.512 & 0.606
\\
& & & & (0.407)& (0.392) & (0.05)%& (0.344) & (0.385)
\\[3pt]
(b) & 0.894 & 0.930 & 0.740 & 0.894 & 0.894 & 0.154 %& 0.540 & 0.638
\\
& & & & (0.536)& (0.530) & (0.064)%& (0.354) & (0.410)
\\[3pt]
(c) & 0.912 & 0.844 & 0.518 & 0.928 & 0.994 & 0.064 %& 0.484 & 0.594
\\
& & & & (0.252)& (0.247) & (0.017)%& (0.320) & (0.360)
\\[3pt]
(d) & 0.892 & 0.878 & 0.622 & 0.880 & 0.890 & 0.098 %& 0.780 & 0.566
\\
&& & & (0.253)& (0.261) & (0.017) %& (0.142) & (0.200)
\\
\hline
\end{tabular*}    \vspace*{-3pt}
\end{table}
%

%s6.2 #&#
%s6.2 ###
\subsection{Comparison with a perturbation based method}
\label{sim-sub1}

In the $p \leq n$ case, \citet{minnier-jasa-11} suggested a
perturbation-based approach for construction of CIs of underlying regression
parameters, including the zero parameters. We compare the performance
of our proposed bootstrap-based method with their approach. We
use $(n=100, p =10)$. The design vectors $\xb_i$ are\vadjust{\goodbreak}
independently selected from a normal population with mean
$\Zero$, unit variances and pairwise covariances
equal to $0.2$. The errors $\varepsilon_i$ are i.i.d. $N(0,\sigma^2)$.
We considered two choices, $\sigma= 1$ and $5$.
The true regression parameter is $\bb= {(2,-2,0.5,
-0.5,0,\ldots,0)}^\prime$. This is very similar to the setup
used in \citet{minnier-jasa-11}. Among the different types of CIs they
proposed, we focus on (i) the usual normal type CI
(which has been modified by a thresholding approach
to handle underlying zero parameters) and denoted by $\operatorname
{CR}^{*\mathrm{N}}$
and (ii) CIs directly based on the quantiles of
the perturbed regression estimates,
denoted by $\operatorname{CR}^{*\mathrm{Q}}$. As suggested in their paper,
we used a BIC-based choice for $\lambda_{2,n}$ for the
simulations; cf. \citet{minnier-jasa-11}.\looseness=-1

%t2 ###
\begin{table}[b]
\caption{Comparison of empirical coverage probabilities and
average lengths (in parentheses) for~$90\%$ CIs for the underlying
parameter $\beta_4(=0.9)$ in cases \textup{(a)} and \textup{(b)}.\break In both cases
$\lambda_{2,n} = 2n^{1/4}$ and in case \textup{(b)}, $\lambda_{1,n} = 0.5
n^{1/2}$}\label{tab2}
\begin{tabular*}{\textwidth}{@{\extracolsep{4in minus 4in}}lcccccc@{}}
\hline
& \multicolumn{3}{c}{\textbf{One-sided}} & \multicolumn{3}{c@{}}{\textbf{Two-sided (with average
lengths)}}\\[-4pt]
& \multicolumn{3}{c}{\hrulefill} & \multicolumn{3}{c@{}}{\hrulefill}\\
\textbf{Case} & $\bolds{\R_n}$ & $\bolds{\breve{\R}_n}$ & \textbf{Oracle} & $\bolds{\R_n}$ & $\bolds{\breve{\R}_n}$ &
\multicolumn{1}{c@{}}{\textbf{Oracle}}\\ %& $|\R_n|$\tmark[a] & $|\breve{\R}_n|$\tmark[a]
\hline
(a) & 0.868 & 0.946 & 0.840 & 0.902 & 0.944 & 0.086 %& 0.512 & 0.606
\\
&& & &  (0.598) & (0.529) & (0.061)%& (0.344) & (0.385)
\\[3pt]
(b) & 0.908 & 0.944 & 0.904 & 0.886 & 0.942 & 0.072 %& 0.540 & 0.638
\\
&& & &  (0.607)& (0.652) & (0.058)\\%& (0.354) & (0.410)
\hline
\end{tabular*}
\end{table}

%t3 ###
\begin{table}
\caption{Comparison of empirical coverage probabilities for
$90\%$ two-sided CIs using the~perturbation based approach by \citet
{minnier-jasa-11}, the oracle and~the~bootstrap based methods. For the
Oracle and Bootstrap methods, the~penalty~parameter is $\lambda_{2,n} =
0.5\cdot n^{1/4}$ and for the perturbation based\break approach the BIC based
choice of $\lambda_{2,n}$ was used}\label{tab3}
\begin{tabular*}{\textwidth}{@{\extracolsep{\fill}}lcccccc@{}}
\hline
& & \multicolumn{2}{c}{\textbf{Perturbation}} &&
\multicolumn{2}{c@{}}{\textbf{Bootstrap}}\\[-6pt]
& & \multicolumn{2}{c}{\hrulefill} &&
\multicolumn{2}{c@{}}{\hrulefill}\\
\textbf{Parameter} & \multicolumn{1}{c}{$\bolds{\sigma}$} & \multicolumn{1}{c}{$\bolds{\operatorname{CR}^{*\mathrm{N}}}$} & \multicolumn{1}{c}{$
\bolds{\operatorname{CR}^{*\mathrm
{Q}}}$} & \textbf{Oracle} & \multicolumn{1}{c}{$\bolds{\R_n}$} & \multicolumn{1}{c@{}}{$\bolds{\breve{\R}_n}$}\\ %& $\R_n^{*}$ & $\breve{
\hline
$\beta_1 = 4$ & 1 & 0.012 & 0.306 & 0.132 & 0.916 & 0.898 %& -- & --
\\
%& && & && 1 & 0.116 & 0.898 & 0.878 & -- & -- \NN
%& && & && 2 & 0.080 & 0.870 & 0.826 & -- & -- \NN[1.2ex]
& 5 & 0.122 & 0.876 & 0.124 & 0.916 & 0.914 %& -- & --
\\[3pt]
%& && & && 1 & 0.122 & 0.282 & 0.920 & -- & -- \NN
%& && & && 2 & 0.130 & 0.284 & 0.924 & -- & -- \ML
$\beta_5 = 0$ & 1 & 1.0\phantom{00} & 1.0\phantom{00} & 0\phantom{000.} & 0.894 & 0.936 %& 0.894 & 0.936
\\
% & && & && 1 & 0 & 0.062 & 0.102 & 0.062 & 0.102 \NN
% & && & && 2 & 0 & 0.000 & 0.002 & 0.000 & 0.002 \NN[1.2ex]
& 5 & 0.288 & 0.902 & 0\phantom{000.} & 0.932 & 0.918 %& 0.388 & 0.912
\\
\hline
\end{tabular*}        \vspace*{-5pt}
\end{table}
% & && & && 1 & 0 & 0.442 & 0.906 & 0.442 & 0.906 \NN
% & && & && 2 & 0 & 0.524 & 0.908 & 0.524 & 0.908 \LL

As shown in Table~\ref{tab3} and somewhat
contrary to the findings of \citet{minnier-jasa-11}, we found
that the $\operatorname{CR}^{*\mathrm{N}}$ based CIs have poor coverage
for both zero and nonzero regression parameters. However, the
$\operatorname
{CR}^{*\mathrm{Q}}$ method performs much better,
particularly when the error variance is high. In comparison,
the bootstrap-based methods are uniformly superior in all cases.
We also noted that compared to the the $\operatorname{CR}^{*\mathrm
{Q}}$ method,
the coverage accuracy of the bootstrap CIs is more sensitive to the
choice of the smoothing parameter
for the zero parameters; see Section~\ref{sim-sub2} below.

%s6.3 #&#
%s6.3 ###
\subsection{Choice of tuning parameter}
\label{sim-sub2}

For penalized regression techniques, the cross validation (CV) has been
a popular method for choosing the tuning parameters, in both low and
high-dimensional cases. We compare the performance of cross validation
(CV) based and theoretical choices of tuning parameters. Based on the
theoretical rates, we use $\lambda_{2,n}=2n^{1/4}$ (for the ALASSO
stage) and in the $p>n$ case, the tuning parameter $\lambda_{1,n}$,
used for the LASSO stage, is set at $\lambda_{1,n}=0.5n^{1/2}$. When
using CV, the initial tuning parameter $\lambda_{1,n}$ is selected by
5-fold CV (only in the $p>n$ case) and kept fixed. Using this fixed
value and again using 5-fold CV, the tuning parameter $\lambda_{2,n}$
for the ALASSO stage is selected. When the underlying true parameter is
zero, an additional theoretical choice of $\lambda_{2,n}=0.25\cdot
{n^{1/4}}$ is used for comparison.

As seen from Table~\ref{tab4}, in case (a) (with $p<n$), using the
CV-based choice of $\lambda_{2,n}$ leads to very good empirical
coverage probabilities for all choices of underlying regression
parameters, including zero parameters. The theoretical choice also
performs comparably for all parameters, except the zero parameter case,
where a smaller value of $\lambda_{2,n}$ performs comparably. The
results in Table~\ref{tab5}, for case (b) (in the $p>n$ setup), show
that there is an overall decrease in the empirical coverage
probabilities for both choices. Unlike the results in case (a) (cf.
Table~\ref{tab4}), the performance is very poor for the zero parameters
irrespective of the method used for selecting the tuning
parameters.\looseness=-1

%t4 ###
\begin{table}
  \caption{Comparison of empirical coverage probabilities for
90\% CIs for different parameters, using CV based and theoretical
choices of $\lambda_{2,n}$ in case \textup{(a)}. The optimal CV based
$\lambda_{2,n}= 0.049\approx0.017\cdot{60}^{1/4}$. For the zero
parameter case an additional (theoretical) choice of $\lambda
_{2,n}=0.25*n^{1/4}$ is compared}\label{tab4}
\begin{tabular*}{\textwidth}{@{\extracolsep{\fill}}lccccccc@{}}
\hline
& & \multicolumn{3}{c}{\textbf{One-sided}} &
\multicolumn{3}{c@{}}{\textbf{Two-sided}}\\[-6pt]
& & \multicolumn{3}{c}{\hrulefill} &
\multicolumn{3}{c@{}}{\hrulefill}\\
\textbf{Parameter} & \textbf{Method} & $\bolds{\R_n}$ & $\bolds{\breve{\R}_n}$ & \textbf{Oracle} & $\bolds{\R_n}$ & $
\bolds{\breve
{\R}_n}$ & \textbf{Oracle}\\ %& $|\R_n|$\tmark[a] & $|\breve{\R}_n|$\tmark[a]
\hline
$\beta_1=4$ & CV & 0.892 &0.894 &0.588 &0.938 &0.890 &0.162 \\
& Th. & 0.894 &0.898 &0.668 &0.922 &0.894 &0.158\\[3pt]
$\beta_4=0.9$ & CV & 0.882 &0.882 &0.566 &0.924 &0.882 &0.156\\
& Th. & 0.872 &0.944 &0.840 &0.940 &0.864 &0.138 \\[3pt]
$\beta_6=0$ & CV & 0.888 &0.886 &0.428 &0.942 &0.902 & 0\phantom{000.} \\
& Th. & 0.004 &0.004 &0.004 & 0\phantom{000.} &0\phantom{000.} &0\phantom{000.}\\
& Th.\tabnoteref{t1} & 0.896 & 0.850 & 0.180 & 0.944 & 0.884 & 0\phantom{000.} \\
\hline
\end{tabular*}
\tabnotetext[a]{t1}{At $\lambda_{2,n}=0.25*n^{1/4}$.}
\end{table}

%t5 ###
\begin{table}[b]
\caption{Comparison of empirical coverage probabilities for
90\% CIs for different parameters, using CV based and theoretical
choices of $\lambda_{1,n}$ and $\lambda_{2,n}$ in case \textup{(b)}.
The
optimal CV based choices were $\lambda_{1,n}=0.124\approx0.016\cdot
{(60)}^{1/2}$ and $\lambda_{2,n} = 0.639\approx0.229
\cdot{(60)}^{1/4}$}\label{tab5}
\begin{tabular*}{\textwidth}{@{\extracolsep{\fill}}lccccccc@{}}
\hline
& & \multicolumn{3}{c}{\textbf{One-sided}} &
\multicolumn{3}{c@{}}{\textbf{Two-sided}}\\[-6pt]
& & \multicolumn{3}{c}{\hrulefill} &
\multicolumn{3}{c@{}}{\hrulefill}\\
\textbf{Parameter} & \textbf{Method} & $\bolds{\R_n}$ & $\bolds{\breve{\R}_n}$ & \textbf{Oracle} & $\bolds{\R_n}$ & $
\bolds{\breve
{\R}_n}$ & \textbf{Oracle}\\ %& $|\R_n|$\tmark[a] & $|\breve{\R}_n|$\tmark[a]
\hline
$\beta_1=4$ & CV & 0.81\phantom{0} &0.838 &0.730 &0.636 &0.506 &0.104 \\
& Th.&0.894 &0.930 &0.740 &0.894 &0.894 &0.154 \\[3pt]
$\beta_4=0.9$ & CV &0.798 &0.854 &0.748 &0.656 &0.488 &0.104 \\
& Th.&0.908 &0.944 &0.904 &0.886 &0.942 &0.072 \\[3pt]
$\beta_6=0$ & CV &0.384 &0.398 &0.194 &0.216 &0.116 &0.00\phantom{0} \\
& Th.&0.016 &0.016 &0.016 &0\phantom{000.} &0\phantom{000.} &0\phantom{000.} \\
& Th.\tabnoteref{t2} & 0.348 & 0.332 & 0.176 & 0.224 & 0.112 & 0\phantom{000.} \\
\hline
\end{tabular*}
\tabnotetext[a]{t2}{At $\lambda_{2,n}=0.25*n^{1/4}$.}
\end{table}

%s6.4 #&#
%s6.4 ###
\subsection{Real data analysis for the low dimensional case}
\label{sim-sub3}

In this section we apply the bootstrap based methods on a prostrate
cancer data-set, available from a clinical study and used in \citet
{tibs-96} [originally available from \citet{stamey-89}]. In this
clinical study, a total of $n=97$ observations were available and the
variable of interest was log(prostrate specific antigen) ({\texttt
{lpsa}}) and eight different predictors ($p=8$) were used to study the
behavior of this quantity. The predictors were log(cancer volume)
({\texttt{lcavol}}), log(prostrate weight) ({\texttt{lweight}}),
{\texttt{age}}, log(benign prostratic hyperplasia amount) ({\texttt
{lbph}}), seminal vesicle invasion ({\texttt{svi}}), log(capsular
penetration) ({\texttt{lcp}}), Gleason score ({\texttt{gleason}}) and
percentage Gleason scores 4 or 5 ({\texttt{pgg45}}). The columns of the
design matrix are centered and scaled to have unit norm. We use the
following theoretical choice for the penalty parameter: $\lambda
_{2,n}=n^{1/4}$. Table~\ref{real-data-ld} shows CIs for estimated
nonzero coefficients. Note that in more than one instance, the
estimated values of $\beta_{j,n}$ fall outside the bootstrap CIs. This
can be explained by considering that the histograms of the bootstrap
replicates which showed that
% For the regression coefficient $\beta_2$, corresponding to
% the covariate \texttt{lweight}, Figure~\ref{fig1} shows the
%distribution of $\R^{*}_n = \sqrt{n}(\beta^{*}_{2,n} -
% \widehat{\beta}_{2,n})/\sigma^{*}_n$ (see Section~\ref{sec4-sub3}
%for definitions)
% and Figure~\ref{fig2} shows the corresponding bootstrap
%distribution for the bias corrected version $\breve{\R}^*_n$.
the distributions of ${\R}^*_n$ and
$\breve{\R}^*_n$ are
heavily skewed and far from the oracle normal distribution.
This is reflected by the endpoints of the
corresponding CIs in Table~\ref{real-data-ld}.\looseness=1

%t6 ###
\begin{table}
\caption{Analysis of prostrate cancer data from \citet
{tibs-96}. The penalty parameter used~is~$\lambda_{2}= n^{1/4}$. ALASSO
estimates and resultant 90\% two-sided CIs\break for estimated nonzero
components are shown}\label{real-data-ld}
\begin{tabular*}{\textwidth}{@{\extracolsep{\fill}}lcccc@{}}
\hline
\multicolumn{1}{@{}l}{\textbf{Predictor} $\bolds{(j)}$} & $\bolds{\widehat{\beta}_{j,n}}$ & $\bolds{\R_n}$ & $\bolds{\breve{\R}_n}$ &
\multicolumn{1}{c}{\textbf{Oracle}}\\ %& Test
\hline
{\texttt{lcavol}} & 0.688 &(0.520, 0.822) & (0.616, 0.944) & (0.636,
0.741) \\
{\texttt{lweight}} & 0.112 &(0.140, 0.235) & (0.162, 0.395) & (0.067,
0.156)\\
{\texttt{svi}} & 0.167 &(0.138, 0.352) & (0.178, 0.487) & (0.115,
0.219)\\
\hline
\end{tabular*}
\tabnotetext[*]{t3}{Obtained from
\texttt{\href{http://www-stat.stanford.edu/\textasciitilde tibs/ElemStatLearn/datasets/prostate.data}{http://www-stat.stanford.edu/\textasciitilde tibs/ElemStatLearn/datasets/}
\href{http://www-stat.stanford.edu/\textasciitilde tibs/ElemStatLearn/datasets/prostate.data}{prostate.data}}.}
\end{table}

% \begin{minipage}[b]{0.45\linewidth}
% \centering
% \includegraphics[width=\textwidth]{hist1.pdf}
% \caption{default}
% \label{figfigure1}
% \end{minipage}
%
% \begin{minipage}[b]{0.45\linewidth}
% \centering
% \includegraphics[width=\textwidth]{hist2.pdf}
% \caption{default}
% \label{figfigure2}
% \end{minipage}
% \end{figure}

%parameter of interest is $\beta_2$, corresponding to the covariate {

%version $\breve{\R}^{*}_n$ when the parameter of interest is $
%Prostrate cancer data.}

%s6.5 #&#
%s6.5 ###
\subsection{Real data analysis for the high-dimensional case}
\label{sim-sub4}

The data, available from a microarray experiment was collected from
\citet{hall-miller-jcgs} and originally used in \citet{segal-03}. The
data consisted of observations from $n=30$ specimens on the Ro1
expression level~($y$), and genetic expression levels $\xb
={(x_1,\ldots,x_p)}^\prime$ for $6319$ genes. The absolute value of the correlation
between $y$ and each covariate $x_i$ was used as an initial screening
tool and only those covariates with absolute correlation value $\geq$0.5 were selected for further study. This resulted in a smaller set of
$p=545$ covariates. The columns of the design matrix were centered and
scaled (by the columnwise standard deviation) and the response vector
$\yb$ was also transformed by centering and scaling. The selected
tuning parameters were $\lambda_1 = 0.5\cdot n^{1/2}$ and $\lambda_2 =
0.5\cdot{n}^{1/4}$. After the initial LASSO step, twenty covariates
are selected and after the ALASSO step only six covariates (genes) were
selected (shown in Table~\ref{real-data-hd}). The residual sum of
squares divided by $(n\mbox{-number of nonzero parameters})$ provides
the following: for the initial LASSO estimate $0.1082$ (equivalent to a
$R^2$ value of $0.888$) and for the ALASSO estimate we obtain $0.092$
(equivalent to $R^2 = 0.904$). This suggests that the extra 14
variables, present in the LASSO estimator provide very little
information about the response. Note that here also the estimated
values of $\beta_{j,n}$'s
often fall outside the bootstrap CIs based on the bias
corrected pivot $\breve{\R}_n$. This suggests that the
true values of the nonzero parameters are probably much
larger in absolute value than suggested by their ALASSO
point estimates.

%t7 ###
\begin{table}
\caption{Analysis of microarray data with $n=30$ and $p=545$
(after initial screening step).\break All six predictors with nonzero ALASSO
coefficients and corresponding 90\% two-sided CIs based on the
bootstrap and oracle methods}\label{real-data-hd}
\begin{tabular*}{\textwidth}{@{\extracolsep{\fill}}lcccc@{}}
\hline
\multicolumn{1}{@{}l}{\tabnoteref{t4}\textbf{Predictor} $\bolds{(j)}$} & $\bolds{\widehat{\beta}_{j,n}}$ & $\bolds{\R_n}$ & $\bolds{\breve{\R}_n}$ &
\multicolumn{1}{c}{\textbf{Oracle}}\\ %& Test
\hline
{\texttt{G709}} & $-0.066$ & $(-0.146, -0.120)$ &
$(-0.490, -0.331)$ & $(-0.127, -0.005)$\\[2pt]
{\texttt{G2272}} & \phantom{$-$}$0.095$ &$(0.087, 0.207)$ & $(0.376, 0.619)$& $(0.010,
0.180)$\\[2pt]
{\texttt{G3655}} & \phantom{$-$}$0.475$ &$(0.250, 0.759)$ & $(0.749, 1.309)$ & $(0.375,
0.575)$\\[2pt]
{\texttt{G4322}} & $-0.021$ &$(-0.047, -0.041)$ &
$(-0.443, -0.432)$ & $(-0.091, 0.048)$\\[2pt]
{\texttt{G5904}} & \phantom{$-$}$0.240$ &$(0.161, 0.507)$ &$(0.495, 0.900)$ & $(0.168,
0.311)$\\[2pt]
{\texttt{G6252}} & \phantom{$-$}$0.112$ &$(0.029, 0.241)$ & $(0.414, 0.687)$ & $(0.030,
0.193)$\\
\hline
\end{tabular*}
\tabnotetext[a]{t4}{Data available from
supplementary material of \citet{hall-miller-jcgs}.}\vspace*{3pt}
\end{table}

%%%%%%%%%%%%%%%%%%%%%%%%%%%%%%%%%%%%%%%%%%%%%%%%%%%%%%%%%%%%%%%%%%%%%%%%%%%%%%%%%%%%%%%%%%%%%%%%%%%%%%%%%%%%%%%%Sec:
%s7 #&#
%s7 ###
\section{Proofs}
\label{sec6}

%s7.1 #&#
%s7.1 ###
\subsection{Notation}
\label{sec6-sub1}

For notational simplicity, we shall set $p_n=p$, $p_{0,n}=p_0$.
Let $\mathbb{Z}_{+} =
\{0,1,\ldots\}$. Let $K, K(\cdot) \in(0,\infty)$ denote
generic constants not depending on
their arguments (if any), but
not on $n$. Also, in the proofs below, let $n_0\geq1$
denotes a generic (large) integer.
For $\alb= (\alpha_1,
\ldots,\alpha_r)\in\mathbb{Z}^r_{+}$, let $|\alb|=\al_1+\cdots
+\al_p$,
$\alb!=\al_1!\cdots\al_r!$ and let $D^{\alb}$ denote the
differential operator\vspace*{-1pt}
$
\frac{\partial^{|\alb|}}{\partial x^{\alpha_1}_1\cdots\partial
x^{\al_r}_r}
$
on $\mathbb{R}^r$, where $r\geq1$ is an integer. Let
$\W_n=n^{-1/2}\sum_{i=1}^n \xb_i'\varepsilon_i$. Partition $\W_n$
as $\W_n = (\W_n^{(1)'},\W_n^{(2)'})'$, where $\W_n^{(1)}$ is
$p_0\times1$.
Also, set
$\W^{(0)}_n = \W_n$, $p^{(0)}=p$, $p^{(1)} = p_0$ and $p^{(2)} =p- p_0$.
Let $\bbf_n=\D^{(1)}_n\C^{-1}_{11,n}\sbf^{(1)}_n\cdot
\lambda_n n^{-1/2}$, $\Up_n=n^{-1}\sum_{i=1}^n \xib^{0}_i{(\xib
^{0}_i)}{}^{\prime}$\vspace*{2pt}
and
$
\breve{\Up}_n=n^{-1}\sum_{i=1}^n (\xib^{0}_i
+\etab^{(0)}_i){(\xib^{0}_i+\etab^{(0)}_i)}^\prime,
$
where\vspace*{1pt} $\xib^{(0)}_i=\D^{(1)}_n\C^{-1}_{11,n}\xb^{(1)}_i$, $\etab^{(0)}_i
= \D^{(1)}_n\C^{-1}_{11,n}\etab_i$ and $\etab_i
= (\xi_{i,1},\ldots,\xi_{i,p_0} )$ with
$
\xi_{i,j} = -\frac{\lambda_n}{n^{1/2}}\cdot\tilde{x}_{i,j}
\cdot\operatorname{sgn}(\beta_{j,n})\gamma{|\beta
_{j,n}|}^{-(\gamma+1)},
1\leq j\leq p_0
$.\vspace*{2pt}
Next note that by conditions (C.2), (C.3)
and (C.6),
\[
\|\bbf_n\|\leq \bigl\|\D^{(1)}\C_{11,n}^{-1/2}
\bigr\|\cdot\bigl\|\C_{11,n}^{-1/2}\bigr\| \cdot\bigl\|\sbf^{(1)}_n
\bigr\| \cdot\frac{\lambda_n}{\sqrt{n}} = O \bigl(n^{-\delta} \bigr).
\]
Let $r_1=\min\{r\geq1\dvtx\|\bbf_n\|^{r+1} = o(n^{-1/2})\}$.
Define the Lebesgue density of the EE for $\T_n$ by
\begin{eqnarray*}
\psi_n(\xb) &= &\phi \bigl(\xb,\sigma^2\breve{
\Upsilon}_n \bigr) \Biggl[1 + \sum_{|\alb|=1}^{r_1}
\bbf^{\alb}_n \chi_{\alb} \bigl(\xb;
\sigma^2\breve{\Upsilon}_n \bigr)
\\
&&\hphantom{\phi \bigl(\xb,\sigma^2\breve{
\Upsilon}_n \bigr) \Biggl[} {}+\frac{\mu_3}{6\sqrt{n}} \sum_{|\alb|=3}\bar{
\xib}^{(0)}_n(\alb)\chi_{\alb} \bigl(\xb;
\sigma^2 \breve{\Upsilon}_n \bigr) \Biggr],\qquad \xb\in
\mathbb{R}^q,
\end{eqnarray*}
where $\bar{\xib}^{(0)}_n(\alb) = n^{-1}
\sum_{i=1}^n { (\xib^{(0)}_i )}^{\alb}$, $\phi(\xb,\Upsilon)$
denotes the density of the\break $N(\Zero,\Upsilon)$
distribution on $\mathbb{R}^q$ and
where $\chi_{\alb}(\xb;\Upsilon)$ is defined by the identity
\[
\chi_{\alb}(\xb;\Upsilon)\phi(\xb;\Upsilon) = {(-D)}^{\alb}
\phi(\xb;\Upsilon),\qquad \alb\in\mathbb{Z}^q_{+}.
\]
Next define
the density of the EE for $\R_n$ by
%%%%%%%%%%%%%%%%%
%%%%%{\footnote{doubt: check the expression}}
%
\begin{eqnarray}
\pi_n(\xb) & =& \phi(\xb,\breve{\Upsilon}_n ) \Biggl[1 +
\sum_{k=1}^{r_1} \frac{1}{k!} \biggl\{
\sum_{|\alb|=k}{ ( -\bbf_n )}^{\alb}
\chi_{\alb}(\xb\dvtx\breve{\Upsilon}_n) \biggr\}
\nonumber\\
&&\hspace*{44pt}{} +\frac{1}{\sqrt{n}}\cdot\frac{\mu_3}{6\sigma^3} \biggl\{ \sum
_{|\alb|=1}\sum_{|\glb|=2} \bigl[\bar{\xib
}^{(0)}_n (\alb+\glb)-3\bar{\xib}^{(0)}_n(
\alb)\bar{\xib}^{(0)}_n(\glb) \bigr]\nonumber\\
&&\hspace*{150pt}{}\times \chi_{\alb+\glb}(
\xb;\breve{\Upsilon}_n)
\nonumber\\
&&\hspace*{174pt}{}-3\sum_{|\alb|=1}\bar{\xib}^{(0)}_n(
\alb) \chi_{\alb}(\xb;\breve{\Upsilon}_n) \biggr\} \Biggr],\nonumber\\
\eqntext{\xb\in\mathbb{R}^q.}
\end{eqnarray}

%s7.2 #&#
%s7.2 ###
\subsection{Auxiliary results}
\label{sec5-sub3}

%le7.1 #&#
%
\begin{lemma}
\label{lem5-1}
Under
%{\footnote{doubt: check the experssions}}
\textup{(C.2)} and \textup{(C.4)}:
\begin{longlist}[(iii)]
\item[(i)] $\pr(\|\W^{(1)}_n\|>K\sqrt{p_0\log{n}} )
= O (p_0 \cdot n^{-(r-2)/2} )$;
\item[(ii)] $\pr({\|\W^{(l)}_n\|}_{\infty}>K\sqrt{\log
{n}} )
= O (p^{(l)}\cdot n^{-(r-2)/2} )$, for $l=0,1,2$;
\item[(iii)] $\pr({\|\sqrt{n}(\bbtn
-\bb_n)\|}_{\infty}>K\sqrt{\log{n}} )=O (p\cdot
n^{-(r-2)/2} )$.
\end{longlist}
\end{lemma}

\begin{pf}
See the supplementary material \citet{CL} (hereafter referred
to as [CL]).
\end{pf}

The key step in the proofs of Theorems~\ref{thm3-1}--\ref{thmhd-1}
is EEs for the ALASSO estimator and its studentized
version which are given below.

%th7.2 #&#
%
\begin{theorem}
\label{thm5-1}
\textup{(a)} If conditions \textup{(C.1)--(C.6)}
hold with $r=4$, then
\[
\sup_{B\in\cC_q} \biggl| \pr(\T_n\in B ) - \int
_{B} \psi_n(\xb) \,d\xb\biggr|  = o
\bigl(n^{-1/2} \bigr).
\]

\begin{longlist}[(b)]
\item[(b)] If conditions \textup{(C.1)}$^{\prime}$--\textup{(C.6)}$^{\prime}$ hold
with $r=6$, then
\[
\sup_{B\in\cC_q}\biggl | \pr(\R_n\in B ) - \int
_{B} \pi_n(\xb) \,d\xb\biggr|  = o
\bigl(n^{-1/2} \bigr).
\]
\end{longlist}
\end{theorem}

\begin{pf}%[Proof of Theorem~\ref{thm5-1}]
See [CL].
\end{pf}

%s7.3 #&#
%s7.3 ###
\subsection{Proof of the main results}
\mbox{}
\begin{pf*}{Proof of Theorem~\ref{thm3-1}}
We only\vspace*{1pt} give an outline of the proof here. For the details of the steps,
see [CL].
Let $\bolds{\La}^{(1)}_n$ be a $p_0\times p_0$
diagonal matrix with $j$th
diagonal entry given by
$\operatorname{sgn}(\beta_{j,n}){|\beta_{j,n}|}^{-(\gamma+1)}$,
$1\leq j\leq p_0$. Then it can be shown that
%
%e7.1 #&#
%
%e7.1 ###
\begin{equation}
n^{-1}\sum_{i=1}^n
\xib^{(0)}_i{\etab^{(0)}_i}^\prime
 = -\frac{\lambda_n\gamma}{n}\D^{(1)}_n \C^{-1}_{11,n}
\bolds{\La}^{(1)}_n\C^{-1}_{11,n}{\D
^{(1)}_n}^{\prime}. \label{p31-1}
\end{equation}
Using
Theorem~\ref{thm5-1}(a),
one gets
%
%e7.2 #&#
%
%e7.2 ###
\begin{eqnarray}
\label{p31-2} \Delta_n & \equiv&\sup_{B\in\cC_q} \biggl|\pr(
\T_n\in B ) -\int_{B}\phi \bigl(\xb;
\sigma^2\Upsilon_n \bigr) \,d\xb\biggr|
\nonumber
\\
%& = &\sup_{B\in\cC_q}\big|\int_B\left[\Psi_n(\xb)-\phi(\xb'\sigma^2
% \right] d\xb\big|+o\left(n^{-1/2}\right)\nonumber\\
& = &\sup_{B\in\cC_q} \biggl|\int
_B \bigl[\phi \bigl(\xb; \sigma^2\breve{
\Upsilon}_n \bigr)-\phi \bigl(\xb;\sigma^2
\Upsilon_n \bigr) \bigr] \,d\xb
\nonumber
\\
& & \hspace*{22pt}{}+\sum_{|\alb|=1}\bbf^{\alb}_n
\int_B \chi_{\alb} \bigl(\xb;\sigma^2
\breve{\Upsilon}_n \bigr) \phi \bigl(\xb;\sigma^2\breve{
\Upsilon}_n \bigr) \,d\xb
\nonumber
\\[-8pt]
\\[-8pt]
\nonumber
& &\hspace*{22pt} {}+ \frac{\mu_3}{6\sqrt{n}}\sum_{|\alb|=3} \bar{
\xib}^{(0)}_n(\alb)\int_B
\chi_{\alb} \bigl(\xb; \sigma^2\breve{\Upsilon}_n
\bigr)\phi \bigl(\xb; \sigma^2\breve{\Upsilon}_n \bigr)
\,d\xb \biggr|
\\
& & {}+o \bigl(n^{-1/2}+\|\bbf_n\| \bigr)
\nonumber
\\
& \equiv&\sup_{B\in\cC_q}\bigl|I_{1,n}(B)+I_{2,n}(B)+I_{3,n}(B)\bigr|
+o \bigl(n^{-1/2}+\|\bbf_n\| \bigr).\nonumber
\end{eqnarray}
Also, by
conditions (C.2)--(C.6),
%
%e7.3 #&#
%
%e7.3 ###
\begin{eqnarray}
\label{p31-3} \|\breve{\Upsilon}_n-\Upsilon_n\| & =&
\Biggl\|2n^{-1}\sum_{i=1}^n
\xib^{(0)}_i{\etab^{(0)}_i}^{\prime}
+ n^{-1}\sum_{i=1}^n
\etab^{(0)}_i{\etab^{(0)}_i}^{\prime}
\Biggr\|
\nonumber
\\[-8pt]
\\[-8pt]
\nonumber
& \leq& K(q,\gamma)\cdot\frac{\lambda_n}{n}\cdot n^{a+b(\gamma+1)}.
\end{eqnarray}
The proof of Theorem~\ref{thm3-1} now follows
from \eqref{p31-1}--\eqref{p31-3}; See [CL].
\end{pf*}

\begin{pf*}{Proof of Theorem~\ref{thm3-1b}}
Since
%{\footnote{doubt:what is $\Gamma_n$}}
$\operatorname{tr}(\bolds{\Gamma}_n)\geq\delta q n^{a+b(\gamma
+1)}$ for
some $\delta\in(0,1)$ and $\bolds{\Gamma}_n$ is $q\times q$, for
each $n\geq1$, there exist a $j_n\in\{1,\ldots,q\}$
such that ${(\bolds{\Gamma}_n)}_{j,j}\geq\delta n^{a+b(\gamma+1)}$.
Write $\cC_{q,n}= \{\{\xb\in\mathbb{R}^q\dvtx x_{j_n}\in(-a,a)\}
\dvtx a\in\mathbb{R} \}$.
Also, let
$
\breve{\tau}^2_n = \sigma^2\cdot{(\breve{\Upsilon}_n)}_{j_n,j_n}
\mbox{ and }
\tau^2_n = \sigma^2\cdot{(\Upsilon_n)}_{j_n,j_n}$.
Then, $I_{k,n}=0$, for all $B\in\cC_{q,n}$ for $k=2,3$,
\eqref{p31-2}
and by \eqref{p31-1}--\eqref{p31-3},
%
%e7.4 #&#
%
%e7.4 ###
\begin{eqnarray}
\label{p31-4} \Delta_n & \geq&\sup_{B\in\cC_n}\bigl|I_{1,n}(B)\bigr|
+ o \bigl(n^{-1/2}+\|\bbf_n\| \bigr)
\nonumber
\\
&=&\sup \biggl\{ \biggl|\int_{-a}^a \bigl[\phi(x,
\breve{\tau}) - \phi(x, \tau) \bigr] \,dx\biggr | \dvtx a\in\mathbb{R} \biggr\}
+o \bigl(n^{-1/2}+\|\bbf_n\| \bigr)
\nonumber
\\[-8pt]
\\[-8pt]
\nonumber
& \geq& K\bigl|\breve{\tau}^2_n -\tau^2_n\bigr|
+o \bigl(n^{-1/2}+\|\bbf_n\| \bigr)
\\
& \geq& K\cdot\delta\gamma\cdot\frac{\lambda_n}{n}\cdot n^{a+b(\gamma+1)} +o
\bigl(n^{-1/2}+\|\bbf_n\| \bigr).\nonumber
\end{eqnarray}
This proves part (b) in the case where
$
n^{-1/2}+\frac{\lambda_n}{\sqrt{n}}\cdot n^{b\gamma}
= O(\lambda_n\cdot n^{-1+a+b(\gamma+1)} )
$. A subsequence argument proves part (b)
when this condition fails. See [CL] for more details.
\end{pf*}

%le7.3 #&#
%
\begin{lemma}
\label{lem5-2}
Suppose that conditions \textup{(C.1)}$^\prime$--\textup{(C.6)}$^\prime$ holds with $r=5$, and let
$n^{-1}\sum_{i=1}^n {\|\C^{-1/2}_{11,n}\xb^{(1)}_i\|}^5=O(1)$.
Then, for any $\delta>0$ and $K\in(0,\infty)$,
there exists $\delta_0\in(0,1)$ such that
\[
\sup \bigl\{\bigl|\widehat{\omega}_n(t_1,t_2)\bigr|
\dvtx\delta^2 \leq t^2_1+t^2_2
\leq n^K \bigr\} = 1-\delta_0+o_p(1),
\]
where
\begin{eqnarray*}
\widehat{\omega}_n(t_1,t_2)&=&\Es\exp{ \bigl(
\iota t_1\varepsilon^{*}_1 + \iota
t_2 { \bigl(\varepsilon^{*}_1
\bigr)}^2 \bigr)},
\\
\omega(t_1,t_2) &=& \E\exp{ \bigl(\iota t_1
\varepsilon_1 +
\iota t_2 {(\varepsilon_1)}^2 \bigr)},\qquad
t_1,t_2\in\mathbb{R}.
\end{eqnarray*}
\end{lemma}

\begin{pf}%[Proof of Lemma~\ref{lem5-2}]
See [CL].
\end{pf}

\begin{pf*}{Proof of Theorem~\ref{thm3-2}}
Restricting attention to a suitable set $A_{3,n}$
with $\pr(A_{3,n})\raw1$ and retracing the steps in
the proof of Theorem~\ref{thm5-1},
one can show (cf. [CL]) that
%
%e7.5 #&#
%
%e7.6 ###
%e7.5 ###
\begin{eqnarray}
\label{p32-2} %
\sup_{B\in\cC_q}
\biggl|\ps \bigl(\T^{*}_n\in B \bigr) -\int_B
\widehat{\psi}_n(\xb) \,d\xb\biggr| & = &o \bigl(n^{-1/2} \bigr);
\nonumber
\\[-8pt]
\\[-8pt]
\nonumber
\sup_{B\in\cC_q} \biggl|\ps \bigl(\R^{*}_n\in B
\bigr) -\int_B\widehat{\pi}_n(\xb) \,d\xb\biggr| & =&
o \bigl(n^{-1/2} \bigr),
\end{eqnarray}
where $\widehat{\psi}_n$ and $\widehat{\pi}_n$ are obtained
from $\psi_n$ and $\pi_n$, respectively, by replacing
$(\sigma^2,\mu_3,\bbf^{\prime}_n)$ by $(\widehat{\sigma}^2_n,
\widehat{\mu}_{3,n},\widehat{\bbf}^{\prime}_n)$, where
\[
\widehat{\sigma}^2_n = \Var_{*} \bigl(
\varepsilon^{*}_1 \bigr),\qquad \widehat{\mu}_{3,n}=
\Es{ \bigl(\varepsilon^{*}_1-\Es\varepsilon^{*}_1
\bigr)}^3,\qquad \widehat{\bbf}_n=\D^{(1)}_n
\C^{-1}_{11,n}\widehat{\sbf}^{(1)}_n,
\]
and the $j$th element of $\widehat{\sbf}^{(1)}_n$ is given by
$\operatorname{sgn}(\widehat{\beta}_{j,n})\lambda_n \cdot
n^{-1/2}\cdot
{|\widehat{\beta}_{j,n}|}^{-\gamma}$, $1\leq j\leq p_0$. For part (a),
we have, for $n\geq n_0$,
\begin{eqnarray*}
&& \pr \Bigl(\sup_{B\in\cC_q} \bigl|\ps \bigl(\T^{*}_n
\in B \bigr) -\pr(\T_n\in B) \bigr|>Kn^{-1/2} \Bigr)
\\
&&\qquad \leq\pr \Bigl( \Bigl\{\sup_{B\in\cC_q} \bigl|\widehat{
\Psi}_n(B) -\Psi_n(B) \bigr|>Kn^{-1/2} \Bigr\}\cap
A_{3,n} \Bigr)+ \pr \bigl(A^c_{3,n} \bigr)
\\
&&\qquad \leq\pr \biggl(\int\bigl|\phi \bigl(\xb;\widehat{\sigma}^2 \breve{
\Upsilon}_n \bigr)-\phi \bigl(\xb;\sigma^2_n
\breve{\Upsilon}_n \bigr) \bigr| \,d\xb>Kn^{-1/2} \biggr)+\pr
\bigl(A^c_{3,n} \bigr)
\\
&&\qquad \leq\pr \bigl( \bigl|\hat{\si}^2_n - \si^2\bigr|
>Kn^{-1/2} \bigr)+ o(1),
\end{eqnarray*}
which can be made arbitrarily small by choosing
$K\in(0,\infty)$ large. Hence, part (a) follows.
The proof of part (b) is similar; see [CL] for more details.~%
\end{pf*}

\begin{pf*}{Proof of Theorem~\ref{thm3-3}}
From the proof of Theorem~\ref{thm5-1}
in [CL], there exists a set $A_{1,n}$
with $P(A_{1,n}^c)= o(n^{-1})$, such that on
$A_{1,n}^c$
and for $n\geq n_0$,
\begin{eqnarray*}
\widehat{I}_n &=& I_n \quad\mbox{and}
\\
\breve{\R}_n & \equiv&\frac{\sqrt{n}\D_n (\bbln-\bb_n )
+\breve{\bbf}_n}{\breve{\sigma}_n}
\\
& =& \biggl[ \biggl\{\D^{(1)}_n\C^{-1}_{11,n}
\W^{(1)}_n -\frac{\lambda_n}{\sqrt{n}}\D^{(1)}_n
\C^{-1}_{11,n} \tilde{\sbf}^{(1)}_n \biggr
\} +\frac{\lambda_n}{\sqrt{n}}\D^{(1)}_n\C^{-1}_{11,n}{
\sbf^{\dagger
}_n}^{(1)} \biggr]\cdot\frac{1}{\breve{\sigma}_n}
\\
& \equiv&\D^{(1)}_n\C^{-1}_{11,n}
\W^{(1)}_n\cdot\frac{1}{\breve{\sigma}_n} + \Q_{3,n}\qquad
\mbox{(say)},
\end{eqnarray*}
where,
$
\Q_{3,n} = \frac{\lambda_n}{\sqrt{n}}\cdot\D^{(1)}_n
\C^{-1}_{11,n} ({\sbf^{\dagger}_n}^{(1)}-\widetilde{\sbf
}^{(1)}_n ),
$
and the $j$th element of ${\sbf^{\dagger}_n}^{(1)}$ is given by
$
s^{\dagger}_{j,n} = \operatorname{sgn}(\widehat{\beta}_{j,n})
{|\widetilde{\beta}_{j,n}|}^{-\gamma}, 1\leq j\leq p_0
$. Note that
\begin{eqnarray*}
& &\pr\bigl(\|\Q_{3,n}\|\neq0 \bigr)
\\
&&\qquad \leq\pr \bigl( \bigl\{{\sbf^{\dagger}_n}^{(1)}\neq
\widetilde{\sbf}^{(1)}_n \bigr\}\cap A_{1,n} \bigr)+
\pr \bigl(A^c_{1,n} \bigr)
\\
& &\qquad\leq\pr \bigl( \bigl\{\operatorname{sgn}(\widehat{\beta}_{j,n}) \neq
\operatorname{sgn}(\beta_{j,n}), \mbox{for some $1\leq j\leq
p_0$ } \bigr\} \cap A_n \bigr)+\pr \bigl(A^c_{1,n}
\bigr)
\\
& &\qquad= 0 + \pr \bigl(A^c_{1,n} \bigr)\qquad \mbox{for $n\geq
n_0$}
\\
&&\qquad = o \bigl(n^{-1} \bigr).
\end{eqnarray*}

Next, using Taylor's expansion,
one can write
\begin{eqnarray*}
\breve{\R}_n & = &\D^{(1)}_n
\C^{-1}_{11,n}\W^{(1)}_n \biggl[
\sigma^{-1} -\frac{1}{2\sigma^3} \bigl(\breve{\sigma}^2_n-
\sigma^2 \bigr) +\frac{3}{4\sigma^5}\frac{{ (\breve{\sigma}^2_n
-\sigma^2 )}^2}{2!} \biggr]+
\Q_{4,n}
\\
& \equiv&\breve{\R}_{1,n}+\Q_{4,n}\qquad \mbox{(say)},
\end{eqnarray*}
where
$
\pr(\|\Q_{4,n}\|>Kn^{-3/2}{(\log{n})}^2 )=o
(n^{-1} )
$.
As a consequence, EEs for $\breve{\R}_n$ and $\breve{\R}_{1,n}$
coincide upto order $n^{-1}$. Now using arguments
in the proof of Theorem~\ref{thm5-1}(b), combined with
the arguments in
\citet{gotze-87} and \citet{lahiri-1994-sankhya},
and then using the transformation technique of \citet{bhat-ghosh-aos-78},
one can show (see [CL] for details) that
%
%e7.6 #&#
%
%e7.7 ###
\begin{equation}
\sup_{B\in\cC_q} \biggl|\pr(\breve{\R}_n\in B ) - \int
_B\pi_{1,n}(\xb) \,d\xb\biggr|  = o
\bigl(n^{-1} \bigr), \label{p33-*}
\end{equation}
where
\[
\pi_{1,n}(\xb)  = \phi(\xb\dvtx\Upsilon_n)
\bigl[1+n^{-1/2}p_{1,n} \bigl(\xb;\sigma^2,
\mu_3 \bigr) +n^{-1}p_{2,n} \bigl(\xb;
\sigma^2,\mu_3,\mu_4 \bigr) \bigr],
\]
with $\mu_4=\E\varepsilon^4_1$ and
where $p_{1,n}(\cdot)$ and
$p_{2,n}(\cdot)$ are polynomials of degree 3 and~6,
respectively, with coefficients that are rational functions of
the respective sets of parameters such that the denominators
depend only on $\sigma^2$ [as in the definition of $\pi_n(\cdot)$].

Next, using Lemma~\ref{lem5-2} % with $k>1$
and similar arguments, one can show that
%
%e7.7 #&#
%
%e7.8 ###
\begin{equation}
\sup_{B\in\cC_q} \biggl|\ps \bigl(\breve{\R}^{*}_n
\in B \bigr) - \int_B\widehat{\pi}_{1,n}(\xb) \,d
\xb \biggr|  = o_p \bigl(n^{-1} \bigr), \label{p33-**}
\end{equation}
where
%$$
%
\[
\widehat{\pi}_{1,n}(x)  = \phi(x;\Upsilon_n) \bigl[1
+n^{-1/2}p_{1,n} \bigl(x;\widehat{\sigma}^2_n,
\widehat{\mu}_{3,n} \bigr) +n^{-1}p_{2,n} \bigl(x;
\widehat{\sigma}^2_n,\widehat{\mu}_{3,n},
\widehat{\mu}_{4,n} \bigr) \bigr],
\]
%
%$$
with $\widehat{\sigma}^2_n = \Es{(\varepsilon^{*}_1)}^2$,
$\widehat{\mu}_{k,n}=\Es{(\varepsilon^{*}_1)}^k$,
$k=3,4$.
% and
%$\widehat{\mu}_{4,n}=\Es{(\varepsilon^{*}_1)}^4$.
Theorem~\ref{thm3-3} now follows from \eqref{p33-*}
and \eqref{p33-**}.
\end{pf*}

\begin{pf*}{Proof of Theorem~\ref{thmhd-1}}
Using the arguments similar to
the proof of Theorem~\ref{thm5-1}, one can show that
%
%e7.8 #&#
%
%e7.9 ###
\begin{equation}
\T_n = \D^{(1)}_n \C^{-1}_{11,n}
\W^{(1)}_n - \bbf_n + \Delta_{1,n}
\equiv\T^{\dagger}_{1,n} + \Delta_{1,n}\qquad \mbox{(say),}
\label{hdpf-1}
\end{equation}
where
%
%e7.9 #&#
%
%e7.10 ###
\begin{equation}
\pr\bigl(\|\Delta_{1,n}\|> K\lambda_n \sqrt{p_0
\log{n}}/n \bigr)  = o \bigl(n^{-1/2} \bigr). \label{hdpf-2}
\end{equation}
Note that by (C.6), $\lambda_n n^{-1}\sqrt{p_0\log{n}} = o(n^{-1/2})$,
when $b=0$.
Now using the arguments in the proof of Theorem~\ref{thm5-1}
(with $\eta^{(0)}_i = 0$ for all $i=1,\ldots,n$),
one can conclude (cf. [CL]) that
%
%e7.10 #&#
%
%e7.11 ###
\begin{equation}
\sup_{B\in\mathcal{C}_q} \biggl|\pr(\R_n\in B) - \int
_B \pi^{\dagger}_n(\xb) \,d\xb\biggr| = o
\bigl(n^{-1/2} \bigr), \label{hdpf-3}
\end{equation}
and that
%
%e7.11 #&#
%
%e7.12 ###
\begin{equation}
\sup_{B\in\mathcal{C}_q} \biggl|\ps \bigl(\R^{*}_n\in B
\bigr) - \int_B { \bigl(\pi^{\dagger}
\bigr)}^{*}(\xb) \,d\xb\biggr|  = o_p \bigl(n^{-1/2}
\bigr), \label{hdpf-4}
\end{equation}
where $\pi^{\dagger}_n(\cdot)$ is defined by setting
$\eta^{(0)}_i=0$ for $1\leq i\leq n$ in the definition of $\pi
_n(\cdot
)$, and
where ${(\pi^{\dagger})}^{*}(\cdot)$ is obtained from $\pi^{\dagger
}(\cdot)$
by replacing $\bbf_n$, $\sigma^2$ and $\mu_3$ with $\widehat{\bbf}_n$,
$\widehat{\sigma}^2$ and $\widehat{\mu}_{3,n}$, as in \eqref{p32-2}.
Using \eqref{hdpf-3} and \eqref{hdpf-4}, one can\vadjust{\goodbreak} conclude that
\[
\sup_{B\in\mathcal{C}_q} \bigl|\pr(\R_n\in B) - \ps \bigl(
\R^{*}_n\in B \bigr) \bigr| = o_p
\bigl(n^{-1/2} \bigr).
\]
The proof for $\breve{\R}_n$ is similar. We omit the routine details to
save space.
\end{pf*}

\section*{Acknowledgments}
We thank three anonymous referees, the
Associate Editor
and the Co-Editor, Professor Tony Cai, for a number of constructive
comments that, in particular, led to the addition of Section \ref
{sec-new-th} on the $p>n$ case and, also the
real data example in Section~\ref{sim-sub4}.

The first author acknowledges the help from the staff, excellent
infrastructure and atmosphere and financial support from the
Statistical and Applied Mathematical Sciences Institute (SAMSI),
Raleigh, NC, and the Department of Statistics at North Carolina State
University, Raleigh, NC, where part of this work was completed.

\begin{supplement}[id=suppA]
\stitle{Supplement to ``Rates of convergence of the Adaptive LASSO estimators
to the Oracle distribution and higher order refinements by the bootstrap''}
\slink[doi]{10.1214/13-AOS1106SUPP} %[doi,text={...}] - jei reikia
%suskaldyti doi
\sdatatype{.pdf}
\sfilename{aos1106\_supp.pdf}
\sdescription{Detailed proofs of all results.}
\end{supplement}

% imsref loaded by akundreckaite, 2013-04-15 12:49:36
%

\printaddresses

\end{document}